% This is in amstex %

\input amstex
\documentstyle{amsppt}
 
 \NoBlackBoxes
\magnification1000
\pagewidth{6.5 true in}
\pageheight{8.5 true in}
\topmatter
\title
Prime factors of dynamical sequences
\endtitle
\author
Xander Faber and Andrew Granville
\endauthor

\abstract  Let $\phi(t)\in \Bbb Q(t)$ have degree $d\geq 2$. For a given rational number $x_0$,
define $x_{n+1}=\phi(x_n)$ for each $n\geq 0$.  If this sequence is not
eventually periodic, and if $\phi$ does not lie in one of two explicitly determined affine conjugacy classes of rational functions, then $x_{n+1}-x_n$ has a primitive prime factor in its numerator for all sufficiently large $n$. The same result holds for the exceptional maps provided that one looks for primitive prime factors in the denominator of $x_{n+1} - x_n$. Hence the result for {\sl each} rational function $\phi$ of degree at least~$2$ implies (a new proof)
that there are infinitely many primes. The question of primitive prime factors of $x_{n+\Delta} - x_n$ is also discussed for~$\Delta$ uniformly bounded.
\endabstract
\thanks{The authors are partially supported by  NSERC, and the first author is supported by an NSF postdoctoral fellowship.}
\endthanks
\address{ \hfill \break
\noindent Xander Faber: \ McGill University, Montr\'eal QC H3A 2K6, Canada.\hfill \break
Email:  xander\@math.mcgill.ca \hfill \break \smallskip
\noindent Andrew Granville: Universit\'e de Montr\'eal,
Montr\'eal QC H3C 3J7, Canada.\hfill \break
Email:  andrew\@DMS.UMontreal.CA}
\endaddress

% \keywords Arithmetic Dynamical System, Primitive Prime Factor, Prime Divisor \endkeywords
% MSC 2010: 37P05 (arithmetic dynamics: polynomial and rational maps), 11A41 (primes), 11B37 (recurrences),  11D59 (Thue/Mahler Equations)

\endtopmatter
\document

\head 1. Introduction \endhead

For a  given sequence of nonzero integers $\{ x_n\}_{n\geq 0}$, a {\sl primitive prime factor}
of $x_n$ is a prime $p_n$ that divides $x_n$ but does not divide any term $x_m$ with $0 \leq m < n$. For example, the nonzero terms of the Fibonacci sequence $a_0=0, a_1=1$ and
$a_{n+2}=a_{n+1}+a_n$ have a primitive prime factor for every $n>12$. (See [3].) 

We call $p_n$ {\sl a super-primitive prime factor} of $x_n$
if $p_n\nmid x_m$ for all $m\ne n$. The Fermat numbers
$F_n=2^{2^n}+1$ are pairwise coprime, and so have a super-primitive prime factor
for every $n\geq 0$.

Given a rational function $\phi(t)$ and a point $x_0\in \Bbb C\cup \{ \infty\}$, we define
$x_{n+1}=\phi(x_n)$ for each $n\geq 0$. 
% If $x_n\in \Bbb Q$ we write $x_n=u_n/v_n$ where $u_n$ and $v_n$ are coprime rational integers.
We have already seen an example, namely the sequence of Fermat numbers. Indeed, note that 
	$$F_{n+1}-2=2^{2^{n+1}}-1=(2^{2^n}+1)(2^{2^n}-1)=F_n(F_n-2),$$
and so if $\phi(t)=t^2-2t+2$ and $F_0=3$, then $F_{n+1}=\phi(F_n)$ for each $n\geq 0$.

Our first goal  was to show that, for any $\phi(t)\in \Bbb Q(t)$, 
the numerator of $x_{n+1}-x_n$ contains a primitive prime factor for all 
sufficiently large $n$ provided $\{ x_n\}_{n\geq 0}$ does not eventually become periodic
(which is equivalent to the statement that the $x_n$ are distinct). 
However this is not always true: Let  $x_0=1$ and
$x_{n+1}=x_n^2/(2x_n+1)$ for $n\geq 0$.  One verifies by induction that
if $F_n$ is the $n$th Fermat number, then
$$
x_n = \frac 1{F_n-2} \ \ \text{and} \ \ x_{n+1}-x_n=-\frac{2^{2^n}}{F_{n+1}-2} \ ,
$$
so that $2$ is the only prime
divisor of the numerator of $x_{n+1}-x_n$ for all $n\geq 0$. Note that the denominator
of $x_n$ has a primitive prime factor for all $n > 0$, which also divides the denominator of
$x_{n+1}$.  This is what we prove in general: 
Define $\Cal F_1$ to be those $\phi(t)\in \Bbb Q(t)$ of the form 
$\sigma^{-1}\circ \psi \circ \sigma$ for
some linear transformation $\sigma(t) = \lambda t+\beta$ with $\lambda\ne 0$, where
	$$
		\psi(t) = \frac{t^2}{t+1} \ \  \text{or} \ \ \frac{t^2}{2t+1} .
	$$

\proclaim{Theorem 1} Suppose that $\phi(t)\in \Bbb Q(t)$ has degree $d\geq 2$,
and that a positive  integer $\Delta$ is given. Let $x_0 \in \Bbb Q$ and define $x_{n+1} = \phi(x_n)$ for each $n \geq 0$. 
If the sequence $\{ x_n\}_{n\geq 0}$ is not eventually periodic, then the numerator of $x_{n+\Delta}-x_n$
has a primitive prime factor for all sufficiently large $n$,
{\rm except} if $\phi\in \Cal F_1$ and $\Delta=1$. (If $x_n$ or $x_{n+\Delta}$ is $\infty$, interpret $x_{n+\Delta} - x_n$ as the ``fraction'' $1/0$.) In the case $\phi \in \Cal F_1$ and $\Delta = 1$, the numerator of $x_{n+1}-x_n$ has the same prime factors for all $n$, and for all sufficiently large $n$ the denominator of $x_n$ has a primitive prime factor, which also divides the 
denominator of $x_{n+1}$.
\endproclaim

The two affine conjugacy classes of rational functions in $\Cal F_1$ are exceptions to the theorem for dynamical reasons. They are characterized by the fact that the point at infinity is fixed with small multiplicity (one or two), and all of the other fixed points are totally ramified. 
% It turns out that totally ramified points don't contribute new prime factors, while the point at infinity will only contribute prime factors to the numerator if it has high multiplicity. 
See Lemma~3 and the proof of Theorem~1.

We call $p_{\Delta, n}$ a {\sl doubly primitive prime factor}  if $p_{\Delta, n}$ divides
the numerator of
$x_{n+\Delta} - x_n$, and if $N\geq n$ and $D\geq \Delta$ whenever 
$p_{\Delta, n}$ divides the numerator of $x_{N+D} - x_N$.
Ingram and Silverman [5, Conj.~20] conjectured that the numerator of
$x_{n+\Delta} - x_n$ has a doubly primitive prime factor for all $\Delta\geq 1$ and $n\geq 0$, 
other than for finitely
many exceptional pairs $(\Delta,n)$. 
Unfortunately their conjecture is false with $\Delta=1$
for any $\phi\in \Cal F_1$, though we believe that an appropriate
modification is true: Let $\Cal B_{\Delta, d}$ be the set of all $\phi(t) \in \Bbb C(t)$ of degree~$d$ such that $\phi$ has no periodic point of exact period $\Delta$. A result of Baker (see \S7 and Appendix~B) shows that $\Cal B_{\Delta, d}$ is non-empty if and only if $(\Delta, d)$ is one of the pairs $(2,2), (2, 3), (2, 4)$ or $(3, 2)$. Define $\Cal F_2$ to be the union of $\Cal B_{2, d}$ for $d = 2, 3, 4$ along with all rational maps $\phi(t)$ of the form $\phi = \sigma^{-1} \circ \psi \circ \sigma$ for some $\sigma(t) = (\alpha t - \beta) / (\gamma t - \delta)$ with $\alpha \delta - \beta \gamma \not= 0$ and $\psi(t) = 1 / t^2$. Define $\Cal F_3 = \Cal B_{3,2}$.  Corollary~2 in \S7 shows that the classes $\Cal F_2$ and $\Cal F_3$ provide further counterexamples to the conjecture of Ingram and Silverman, and we believe that should be all of them.
%We define $p_{\Delta,n}$ to be a {\sl doubly primitive prime factor} of $x_{n+\Delta} - x_n$ if $p_{\Delta,n}$ divides the numerator $x_{n+\Delta} - x_n$, but does not divide the numerator of $x_m-x_\ell$ for any $0\leq \ell < m\leq n+\Delta$ other than when $m=n+\Delta$ and $\ell=n$.

\proclaim{Conjecture} Suppose that $\phi(t)\in \Bbb Q(t)$ has degree $d\geq 2$.
Let $x_0 \in \Bbb Q$ and define $x_{n+1} = \phi(x_n)$ for each $n \geq 0$, and suppose
that  the sequence $\{ x_n\}_{n\geq 0}$ is not eventually periodic. The numerator of $x_{n+\Delta}-x_n$
has a doubly primitive prime factor for  all $n \geq 0$ and $\Delta \geq 1$,
except for those pairs with $\Delta=1, 2$ or $3$ when $\phi\in \Cal F_1, \Cal F_2$ or $\Cal F_3$, respectively, as well as for finitely many other exceptional pairs $(\Delta, n)$. 
\endproclaim

In fact we can prove a strengthening of Theorem 1, which implies that if the above conjecture is false, then there must be exceptional pairs $(\Delta, n)$ with $\Delta$ arbitrarily large. 

\proclaim{Theorem 2} Suppose that $\phi(t)\in \Bbb Q(t)$ has degree $d\geq 2$.
Let $x_0 \in \Bbb Q$ and define $x_{n+1} = \phi(x_n)$ for each $n \geq 0$, and suppose
that  the sequence $\{ x_n\}_{n\geq 0}$ is not eventually periodic. For any given 
$M\geq 1$, the numerator of $x_{n+\Delta}-x_n$ has a doubly primitive prime factor for 
all $n \geq 0$ and $M \geq  \Delta \geq 1$, except  for
those pairs with $\Delta = 1, 2$ or $3$ when $\phi \in \Cal F_1, \Cal F_2$ or $\Cal F_3$, respectively,
as well as for finitely many other exceptional pairs $(\Delta, n)$.
 \endproclaim

Upon iterating the relation $F_{n+1} - 2 = F_n(F_n - 2)$, we see that
$$
F_{n+1}-F_n=2^{2^{n+1}}- 2^{2^n} = 2^{2^n} (F_n-2)=
		2^{2^n} F_{n-1}(F_{n-1}-2)= \cdots = 2^{2^n}F_{n-1}F_{n-2}\cdots F_1F_0.
$$
Hence  $F_{n+1}-F_n$ has the same primitive prime factor $p_{n-1}$ as $F_{n-1}$, but
there can be no super-primitive prime factor since if $p_{n-1}$ divides $F_{n-1}$,
then $p_{n-1}$ divides $F_{N+1}-F_N$ for all $N\geq n$.
On the other hand we saw that all prime factors of $F_n$ are super-primitive,
and this does generalize as we see in the following result
(from which  Theorem 1 and 2 are deduced). For this statement, a point $x_0$ is called
{\sl preperiodic} if and only if the sequence $\{ x_n\}_{n\geq 0}$ is eventually periodic.  

\proclaim{Proposition 1} Let $K$ be a number field.
Suppose that $\phi(t)\in K(t)$ has degree $d\geq 2$,
and that  $0$ is a preperiodic point, but not periodic.
If the sequence of $K$-rationals  $\{ x_n\}_{n\geq 0}$ is not
eventually periodic, then the numerator of $x_n$ has a super-primitive prime (ideal)
factor  $P_n$ for all sufficiently large $n$.\footnote{See section~2 for an explanation of the ``numerator'' of an element of $K$.}
\endproclaim

Although we discovered Proposition~1 independently, we later learned that it appears as a special case of Theorem~7 of [5]. Our general strategy is virtually identical to that of [5], but we have simplified the main Diophantine step in the argument. (Our proof avoids the use of Roth's theorem, and instead proceeds by solving a certain Thue/Mahler equation. We discuss this further at the end of the section.)

We will prove a result analogous to Proposition~1 when $0$ is a periodic point in section~4, though in this case
one does not find super-primitive prime factors. Indeed, if $0$ has period $q$ and
if $P$ divides the numerator of $x_n$, then $P$ divides the numerator of $x_{n+kq}$
for all $k\geq 0$, other than for finitely many exceptional primes~$P$.

The proof of Proposition 1 is based on the following sketch of the
special case $\phi(t)=t^2-2t+2$ and   $x_0\in \Bbb Z$
(which includes another proof for
the special case $x_0=F_0=3$, so that $x_n=F_n$ for all $n\geq 0$).
For this $\phi$ we see that  $0$ is preperiodic but not periodic: $\phi(0)=\phi(2)=2$.  Now
$x_{n+1}=\phi(x_n)\equiv \phi(0) = 2\pmod {x_n}$, and for $m>n$ we then have, by
induction, that $x_{m+1}=\phi(x_m)\equiv \phi(2) = 2\pmod {x_n}$.
Hence if $m>n$, then $(x_m,x_n)=(2,x_n)$ which divides $2$, and so
any odd prime factor of $x_n$ is super-primitive.  If $x_n$ does not have an odd prime factor, then $x_n=\pm 2^k$, and there are only finitely many such $n$ as there are
only finitely many integers $r$ for which  $\phi(r)=\pm 2^k$, and no two $x_n$ can
equal the same value of $r$, else the sequence is eventually periodic.

The deduction of Theorem 1 is based on the following sketch of
the special case $\phi(t)=t^2+3t+1$ and $\Delta=1$, and $x_0 \in \Bbb{Z}$.
Here  $\phi(-2)=\phi(-1)=-1$, so that $-2$ is preperiodic but not periodic.
To be able to apply Proposition 1 we make a linear change of variables
and consider $\psi(t)=\phi(t-2)+2=t^2-t+1$ so that
$\psi(0)=\psi(1)=1$; that is, $0$ is preperiodic but not periodic.
We see that if $y_0=x_0+2$ and $y_{n+1}=\psi(y_n)$, then $y_n=x_n+2$ by induction. So Proposition~1 shows
$y_n=x_n+2$ has a super-primitive prime factor $p_n$ for all $n$ sufficiently large. 
But then $p_n$ divides $x_{n+2} - x_{n+1}$ since
	$$
		x_{n+2} = \phi(\phi(x_n)) \equiv \phi(\phi(-2)) = \phi(-1)= \phi(-2) \equiv \phi(x_n) = x_{n+1} \pmod{p_n},
	$$
while $p_n$ cannot divide $x_{k+1} - x_k$ for any $k \leq n$, else
	$$
		-1 = \phi(-2) \equiv \phi(x_n) = \underbrace{\phi \circ \cdots \circ \phi}_{n-k \text{ times}}(x_{k+1})
			\equiv \underbrace{\phi \circ \cdots \circ \phi}_{n-k \text{ times}}(x_k)  = x_n \equiv -2 
			\pmod{p_n}.
	$$
Hence any super-primitive prime factor $p_n$ of $y_n=x_n+2$ is also a primitive prime factor
of $x_{n+2} - x_{n+1}$.
\iffalse
Now $-2$ is a root of $\phi(t)+1$. If $\phi(t)=-1$ then $\phi(\phi(t))=\phi(-1)=-1=\phi(t)$; and so $t+2$ divides $\phi(t)+1$, which  divides $\phi(\phi(t))-\phi(t)$. Substituting in $t=x_n$ we deduce that $x_n+2$ divides $\phi(\phi(x_n))-\phi(x_n)=x_{n+2}-x_{n+1}$, which implies that a super-primitive prime factor $p_n$ of $x_n+2$ divides $x_{n+2}-x_{n+1}$.  It remains to show that the  $k+1$ is the smallest $k$ for which $p_n$  divides $x_{k+1}-x_{k}$, which can be established using the fact that $\phi(\phi(t))-\phi(t)=(t+1)^2(t+2)^2$
\fi

Having a super-primitive prime factor is, by definition, more rare than
having just a primitive prime factor. At first sight it might seem surprising that one
can prove that esoteric dynamical recurrences have super-primitive prime factors
whereas
second-order linear recurrences (like the Fibonacci numbers) do not. However
the numerator and denominator of the $n$th term of a degree-$d$ dynamical recurrence
grow like $C^{d^n}$, far faster than the $C^n$ of linear recurrences, so we might expect
each new term to have a much better chance of having a prime factor that we have not seen
before.  One approach to proving this is simply based on size, the approach used
for second-order linear recurrences, and so one might believe it should work even
more easily here
--- this is the approach, for instance, of [5].  Our approach uses simple considerations
to imply that a new term in the sequence (or a suitable factor of that term)
can only include  ``old'' prime factors from a finite set, and then
we use the Thue/Mahler theorem to show that this can happen only finitely often. One further
upshot of our method is that it can be made effective; i.e., in principle one could give a bound on the size
of the set of exceptions in our main theorems. Essentially this amounts to applying Baker's method to obtain an effective form of the Thue/Mahler theorem; see the discussion in \S4 for a few more details. 

The remainder of this article is laid out as follows. In the next section, we give a description of the notation used in the paper. In \S3 we give a lower bound on the number of distinct zeros a rational function can have outside of certain exceptional scenarios; this tool will be used in \S4 to prove Proposition~1 and the analogous result for periodic points. As an application we deduce a new unified proof that there are infinitely many primes congruent to 1 modulo a fixed odd prime power. In \S5, \S7 and \S8 we study properties of fixed points and preperiodic points in order to determine when one can change coordinates and apply Proposition~1, and then we use this analysis in \S6 and \S9 to deduce Theorems~1 and~2. Appendix~A contains a number of results from complex dynamics that are used throughout the paper, and Appendix~B recalls the classification of the exceptional rational maps arising from Baker's Theorem, as stated in \S7.

\head 2. Notation \endhead

Suppose that $K$ is a number field, with ring of integers $R$. If $S$ is a finite set of nonzero prime ideals of $R$, write $R_S$ for the ring of $S$-integers --- i.e., the set of all elements $a / b \in K$ where $a,b \in R$ and the ideal $(b)$ is divisible only by primes in $S$. One may always enlarge a given set of primes $S$ so that $R_S$ is a principal ideal domain [2, Prop.~5.3.6]. In that case, any $\alpha \in K$ can be written $\alpha = a/b$ with $a, b \in R_S$ and $(a,b) = (1)$; i.e., $a$ and $b$ share no common prime ideal factor in $R_S$. Moreover, the greatest common divisor of any two elements $a,b \in R_S$, denoted gcd$(a,b)$, is defined to be a generator of the ideal $(a,b)$. It is well defined up a to a unit in $R_S$. 

	We will say that a prime ideal $P$ of $K$ {\it divides the numerator} (resp. {\it denominator}) of an element $\alpha \in K$ to mean that upon writing the fractional ideal $(\alpha)$ as $\frak a / \frak b$ with $\frak a$ and $\frak b$ coprime integral ideals, the ideal $P$ divides $\frak a$ (resp. $\frak b$). When $K = \Bbb Q$ and $P$ is a rational prime number, this agrees with standard usage.

For any field $k$, we identify the projective space $\Bbb P^1(k)$ with $k \cup \{\infty\}$. Returning to the number field $K$ with ring of integers $R$, fix a nonzero prime ideal $P \subset R$. Write $R_P$ for the localization of $R$ at $P$; the ring $R_P$ is the subset of all elements $a/b \in K$ such that $a,b \in R$ and $b \not\in P$. The ring $R_P$ has a unique maximal ideal $PR_P = \{ a/b \in K: a \in P, b \not\in P\}$, and we can identify the residue fields $R/PR \cong R_P / PR_P$. There is a canonical reduction map $\Bbb P^1(K) \to \Bbb P^1(R_P / PR_P)$ given by sending $\alpha \in R_P$ to its image in the quotient $R_P / PR_P$, and by sending $\alpha \in \Bbb{P}^1(K) \smallsetminus R_P$ to $\infty$. We extend the notion of congruences modulo $P$ to $\Bbb P^1(K)$ by saying that $\alpha \equiv \beta \pmod P$ if and only if $\alpha$ and $\beta$ have the same canonical reduction in $\Bbb P^1(R_P / PR_P)$. This gives the usual notion of congruence when restricted to $R_P$; i.e., $\alpha - \beta \in PR_P$ if and only if $\alpha$ and $\beta$ have the same canonical reduction in $\Bbb P^1(R_P / PR_P)$.

For a given $\phi(t)\in K(t)$, fix polynomials $f(t), g(t) \in R[t]$ with no common root
in $\overline{K}$ such that $\phi(t) = f(t) / g(t)$.
%Similarly, define $\phi^{(r)}(t)=f_r(t)/g_r(t)$ for each $r\geq 1$, where
% $f_r(t), g_r(t) \in R[t]$ have no common roots in $\overline{K}$.
Let $d=\deg \phi := \max \{ \deg f, \deg g\}$; then
$F(x,y)=y^d f(x/y),\ G(x,y)=y^d g(x/y)\in R[x, y]$ are homogenous of degree $d$
with no common linear factor over $\overline{K}$. For $F_0(x,y)=x, \ G_0(x,y)=y$, we define
$F_{r+1}(x,y)=F(F_r(x,y), G_r(x,y))$ and
$G_{r+1}(x,y)=G(F_r(x,y), G_r(x,y))$ for all $r\geq 0$.
The polynomials $F_{r}(x,y)$ and $G_{r}(x,y)$ have no common factor in $\overline{K}[x,y]$.

%We will make great use of the fact that if a prime ideal $P$ does not divided both $a, b \in R$, then it doesn't divide both $F_r(a,b)$ and $G_r(a,b)$.

	Throughout we will use the notation $\{x_n\}_{n \geq 0}$ to denote a sequence of elements of $K$ obtained by choosing $x_0 \in K$ and setting $x_{n+1} = \phi(x_n)$ for $n \geq 0$. We will also write $\phi^{(n)}$ for the $n$-fold composition of $\phi$ with itself, so that $x_n = \phi^{(n)}(x_0)$.

\head 3.  Rational functions with many distinct zeros  \endhead

Define  $\Cal T$ to be the set of rational functions $\phi(t)\in \Bbb C(t)$
of degree $d\geq 2$ of one of the following forms:
\smallskip
(i)\  $\phi(t)=t^d/g(t)$ where $g(t)$ is a polynomial of degree $\leq d$;

(ii)\  $\phi(t)=c/t^d$ for some constant $c\ne 0$;
or

(iii)\ $\displaystyle \phi(t)= \frac{\alpha (t-\alpha)^d }{ (t-\alpha)^d-c t^d}$ for some constants
$\alpha, c\ne 0$.
\smallskip

\noindent Note that $0$ is periodic in a period of length one in (i), and in a period of length two
in (ii) and (iii) (where $\phi(0)=\infty$ and $\alpha$, respectively). Geometrically
speaking, $\Cal T$ consists of all rational functions of degree $d$ such that $\phi^{(2)}(t)$ has a totally ramified fixed point at $0$.

The main reason for defining $\Cal T$ is seen in the following lemma:

\proclaim{Lemma 1}  Suppose that $\phi(t)\in \Bbb C(t) \smallsetminus \Cal T$
has degree $d\geq 2$. If $r\geq 4$, then $F_r(x,y)$ has at least three non-proportional linear factors.
\endproclaim

\demo{Proof}   Proposition~3.44 of [12] states that if $\psi(t)\in \Bbb C(t)$ has degree $d\geq 2$, and if
$\psi^{(2)}(t)$ is not a polynomial, then $\psi^{-4}(\infty)$ contains at least three elements.
Let $\psi(t)=1/\phi(1/t)$ so that
$\psi^{(2)}(t) \in \Bbb C[t]$ if and only if $\phi^{(2)}(t)$ is  of the form $t^D/g_2(t)$ where $g_2(t)$ is a non-zero polynomial of degree $\leq D=d^2$. The result follows by showing that this occurs  if and only if $\phi(t)\in \Cal T$:
One easily confirms that, for each $\phi\in \Cal T$, one has
$\phi^{(2)}(t)=t^D/g_2(t)$ for some polynomial $g_2(t)$ of degree $\leq D$. On the
other hand, any  $\phi^{(2)}(t)$ of this form is totally ramified over $0$.
Let $\beta=\phi(0)$ so that $0=\phi(\beta)$.  As ramification indices are multiplicative, we deduce that $\phi$ is totally ramified at $0$ and at $\beta$. An easy calculation then
confirms that the cases $\beta = 0, \infty$, or $\alpha$ ($\ne 0$ or $\infty$)
correspond  to the three cases   in $\Cal T$. \qed
\enddemo

We deduce the following from Lemma 1:

\proclaim{Corollary  1}  Suppose that $\phi(t)\in \Bbb C(t) \smallsetminus \Cal T$ has degree $d\geq 2$. If $r\geq 4$, then $F_r(x,y)$ has at least $d^{r-4}+2$ non-proportional linear factors.
\endproclaim

To prove Corollary~1 we use the $abc$-theorem for polynomials (see, e.g., [4, Thm.~F.3.6]). As we will have use for it again later in the paper, we recall the statement:

\proclaim{$\bold{abc}$-Theorem for Polynomials} If $a, b, c \in \Bbb C[x,y]$ are homogeneous forms of degree $d \geq 1$ with no common linear factor such that $a(x,y) + b(x,y) + c(x,y) = 0$, then the number of non-proportional linear factors of $abc$ is at least $d+2$. 
\endproclaim

\demo{Proof of Corollary~1} We have $F_r(x,y)=F_4(X,Y)$ where $X=F_{r-4}(x,y), Y=G_{r-4}(x,y)$
(which have degree $d^{r-4}$). By Lemma 1, $F_4(X,Y)$ has at least three non-proportional
linear factors in $X,Y$ which must themselves satisfy a linear equation with constant coefficients.
Indeed, if the three linear factors are $X-\alpha Y, X-\beta Y$, and $X-\gamma Y$,
then
$$(\beta-\gamma)(X-\alpha Y)+(\gamma-\alpha)(X-\beta Y)+(\alpha-\beta)(X-\gamma Y)=0.$$
If the three linear factors are instead $X - \alpha Y, X - \beta Y$, and $Y$, then we use
$$ (X - \alpha Y) - (X - \beta Y) + (\alpha - \beta) Y = 0.$$
The $abc$-theorem for polynomials then implies that there are at least $d^{r-4}+2$ coprime linear factors of $(X-\alpha Y) (X-\beta Y) (X-\gamma Y)$, and hence of $F_4(X,Y)=F_r(x,y)$. \qed
\enddemo

%\proclaim{Lemma 2} For any rational function $\phi(t) \in \Bbb C (t)$ and any root $\alpha$ of $H_0(x)$, there
%exists a unique root $\beta$ of $H_0(x)$ with $\phi(\beta)=\alpha$.
%\endproclaim

%\demo{Proof} Let $\beta=\phi^{(\Delta-1)}(\alpha)$ so that
%$\phi(\beta)=\phi\circ \phi^{(\Delta-1)}(\alpha)=\phi^{(\Delta)}(\alpha)=\alpha$
%as $H_0(\alpha)=0$. On the other hand $\beta$ is unique because if
%$\phi(\beta)=\alpha$ then $\phi^{(\Delta-1)} (\alpha)=\phi^{(\Delta)}(\beta)=\beta$.
%\enddemo

%\proclaim{Lemma 3} For any $k\geq 1$ and $\alpha_k$ satisfying
%$H_{k}(\alpha_{k})=0$, we have that $\phi(x+\alpha_{k})-\alpha_{k}\not\in \Cal T$.
%\endproclaim

%\demo{Proof}   If $\phi(x+\alpha_{k})-\alpha_{k} \in \Cal T$ then
%$0$ is a periodic point for $\phi(x+\alpha_{k})-\alpha_{k}$, and so
%$\alpha_k$ is a periodic point for $\phi(x)$. However this contradicts the
%definition of $\alpha_k$.
%\enddemo

\head 4. Preperiodic points and super-primitivity \endhead

We use the following result from Diophantine approximation:

\proclaim{The Thue/Mahler Theorem} Suppose that $F(x,y)\in
K[x,y]$ is homogenous and has at least three non-proportional linear factors
(over $\overline{K}$). Let $S$ be any finite set of primes
of $K$.  There are only finitely many $m/n\in K \smallsetminus \{0\}$ such that all prime
factors of $F(m,n)$ belong to the set $S$.
\endproclaim

See [2, Thm.~5.3.2] for a proof.\footnote{Although we will not need it at present, the Thue/Mahler Theorem can be made effective. One reduces its proof to the solution of a unit equation, and unit equations are effectively solvable by Baker's method. See [2, \S5.4] for a discussion. } With this tool in hand we can complete the proof of Proposition 1:
\medskip

\demo{Proof of Proposition 1} Write $y_0=0$ and $y_{k+1}=\phi(y_k)$
for each $k\geq 0$. Let $r_0 = 0$, $s_0=1$, and choose $r_k, s_k \in
R$ --- the ring of integers of $K$ --- so that $y_k = r_k / s_k$ for each 
$k \geq 0$. If $y_k = \infty$, let $r_k = 1, s_k = 0$.
As $0$ is preperiodic, we may assume there are only finitely many elements
$r_k$ and $s_k$.  Let $S$ be the set of prime ideals that either
divide $r_k$ for some $k \geq 1$ or that divide
Resultant$(F,G)$. Note that $S$ is finite since $r_j\ne 0$ for all $j
\geq 1$, and $f$ and $g$ have no common root in $\bar K$. We may
also enlarge the set $S$ so that the ring of $S$-integers $R_S$ is a principal ideal domain.

Let $P$ be a prime ideal that is not in $S$, and divides the rational
prime $p$, so that $R/PR \cong \Bbb F_q$ for $q$ some power of $p$.
Then $\Bbb P^1(R/PR)\cong \Bbb P^1(\Bbb F_q) = \Bbb F_q \cup \{\infty\}$. We identify $P$ with
the prime ideal $PR_S$ in $R_S$. There is a canonical isomorphism
$R/PR \cong R_S / PR_S$, and hence also $\Bbb P^1(R/PR)\cong \Bbb P^1(R_S / PR_S)$.

At most one term of the sequence $\{x_n\}_{n \geq 0}$ is equal to $\infty$, so we may 
assume that $n$ is large enough that $x_n \not= \infty$. 
 Write each $x_n=u_n/v_n\in K$ with $u_n,v_n\in R_S$ and $(u_n,v_n)=(1)$.
%We will now prove, by induction, that if the prime $P$
%divides $u_n$, then $x_{n+k} \equiv y_k \pmod P$ for all $k \geq 0$. 
% For $k=0$ our claim simply restates that $P$ divides $u_n$
% , since $s_0=1$ and $P \nmid v_n$.  
Suppose $P \mid u_n$.
Observe that since $P$ does not divide Resultant$(F, G)$, for any $k \geq 0$ we find that
$$
	x_{n+k} = \phi^{(k)}(x_n) \equiv \phi^{(k)}(0) = y_k \pmod P.
$$
% let us begin by noting that if $(a,b)=(1)$ then $P$ does not divide gcd$(F(a,b),
%G(a,b))$ (else $P$ divides Resultant$(F,G)$, and so $P\in S$, a
%contradiction), and so if $\phi(a/b)=A/B$ then $(F(a,b):G(a,b))=(A:B)$
%in $\Bbb P^1(R/PR)$.  Now suppose $(u_{n+k}:v_{n+k})=(r_k:s_k)$, so
%there exists an element $\tau\in R_P$, not divisible by $P$,
%with $u_{n+k}\equiv \tau r_k \pmod P$ and $v_{n+k}\equiv \tau s_k
%\pmod P$, and hence $F(u_{n+k},v_{n+k})\equiv \tau^d F(r_k,s_k) \pmod
%P$ and $G(u_{n+k},v_{n+k})\equiv \tau^d G(r_k,s_k) \pmod P$. Therefore
%if $(u_{n+k}:v_{n+k})=(r_k:s_k)$, then in $\Bbb P^1(R/PR)$ we have
%$
%(u_{n+k+1}:v_{n+k+1})=(F(u_{n+k},v_{n+k}): G(u_{n+k},v_{n+k}))
%=(\tau^d F(r_k,s_k):\tau^d G(r_k,s_k))=(r_{k+1}:s_{k+1})$.

We deduce that if a prime ideal $P$ is not in $S$, then $P$ divides at
most one $u_n$.  For if  $P$ divides $u_m$ and $u_n$ with $n>m$,
then $P|r_k$ with $k = n-m$ by the previous paragraph. Hence $P\in S$,
which is a contradiction.

Let $\Cal N$ be the set of integers $n\geq 0$ such that all prime factors of
$u_n$ are in the set $S$. Now $F_4(u_{n-4},v_{n-4})$ equals $u_n
\times$ gcd$(F_4(u_{n-4},v_{n-4}), G_4(u_{n-4},v_{n-4}))$ so that if $n\in \Cal N$ then
all prime factors of $F_4(u_{n-4},v_{n-4})$ are in the set $S$.  Note that
$\phi(t)\not\in \Cal T$ as $0$ is not a periodic point.  Hence, by
Lemma 1, $F_4$ has at least three non-proportional linear factors, and therefore, by the
Thue/Mahler theorem, there are only finitely many 
$u_{n-4}/v_{n-4}\in K$ such that all prime factors of $F_4(u_{n-4},v_{n-4})$
belong to the set $S$.  We deduce that there are only finitely many  $n\in \Cal N$,
since $\{ x_n\}_{n\geq 0}$ is not eventually periodic.

Finally, if $m\not\in \Cal N$, then $u_m$ has a prime ideal factor $P_m$ not in $S$,
and we have seen that $P_m$ cannot divide $u_n$ for any $n\ne m$. Moreover,
$P_m$ cannot divide $v_m$ since $(u_m, v_m)=(1)$. \qed
\enddemo

\remark{Remark 0} The proof of Proposition~1 shows that we may relax the
hypothesis ``0 is preperiodic for $\phi$'' with the possible cost of losing some primitive prime factors.  
More precisely:
Let $\phi  \in K(t) \smallsetminus \Cal T$ be a rational function of degree at least 2, 
and let $x_0$ be a rational number such that the sequence $\{x_n\}_{n \geq 0}$ is not eventually
periodic. Then there exist infinitely many $n$ such that the numerator of $x_n$ admits a primitive
prime ideal factor. In particular, the set $\{P : x_n \equiv 0 \pmod P  \text{ for some }$n$ \}$ is infinite. 
\endremark

\smallskip

Finally we prove a version of Proposition~1 in the case that $0$ is periodic. 
It is also a special case of Theorem~7 of [5], again with the added benefit that 
our method can be made effective. 

\proclaim{Proposition $1'$} Let $K$ be a number field.
Suppose that $\phi(t)\in K(t) \smallsetminus \Cal T$ has degree $d\geq 2$,
and that  $0$ is a periodic point.
If the sequence of $K$-rationals  $\{ x_n\}_{n\geq 0}$ is not
eventually periodic, then the numerator of $x_n$ has a primitive prime factor
$P_n$ for all sufficiently large $n$.
\endproclaim

\remark{Remark 1} It is not hard to see that if $\phi \in \Cal T$, then the conclusion of
Proposition~$1'$ does not follow. Indeed a prime  dividing $x_n$ for some $n\geq 0$,
divides $x_0$ in $\Cal T$(i), 
% and $x_1$ in $\Cal T$(ii), 
and $x_\delta$ in $\Cal T$(ii) and $\Cal T$(iii), 
where $\delta$ is the least non-negative residue of $n\pmod 2$.
\endremark

\remark{Remark 2}  In the proof of Proposition~$1'$ we actually construct a sequence
$\{ u_n^*\}_{n\geq 0}$ of algebraic integers of $K$,
where the numerators $u_n$ of $x_n$ are a  product of powers of $u_m^*$ with $m\leq n$;
and we show that $u_n^*$ has a super-primitive prime factor
$P_n$ for all sufficiently large $n$.  Proposition~$1'$ then follows.
\endremark

\remark{Remark 3} Ingram and Silverman [5] conjecture that if 
$\phi(t)\in K(t)$ has degree $d\geq 2$,
and if $0$ is {\sl not} a preperiodic point,
then for the sequence of $K$-rationals  $\{ x_n\}_{n\geq 0}$ 
starting with $x_0=0$,  the numerator of $x_n$ has a primitive prime factor
for all sufficiently large $n$. This does not seem approachable using our methods.
\endremark

\demo{Proof of Proposition~$1'$}
This is largely based on the above proof when $0$ is preperiodic and not periodic,
but has some additional complications. Suppose that $0$ has period~$q$. Then
$F_q(x,y)=x^b F^*_q(x,y)$ for some integer $b\geq 1$ and $F^*_q(x,y)\in R[x,y]$
where $F^*_q(0,1)\ne 0$. Therefore if $n>q$ then
$$
	F_n(x,y)=F_q(F_{n-q}(x,y), G_{n-q}(x,y))=\Big(F_{n-q}(x,y)\Big)^b F^*_q(F_{n-q}(x,y),G_{n-q}(x,y)).
$$
So we define $F_n^*(x,y)=F_n(x,y)/F_{n-q}(x,y)^b\in R[x, y]$, and we let $F_n^*(x,y)=F_n(x,y)$
if $n<q$. Note that
$$
	F_n(x,y)=\prod_{0\leq j \leq [n/q]} F_{n-jq}^*(x,y)^{b^j}.
$$

Write $y_0=0$ and $y_{k+1}=\phi(y_k)$ for each $k\geq 0$. Let $S$ be the set of primes of $K$ that either divide the numerator of $y_k$ for some $y_k \not= 0$, or that divide Resultant$(F,G)$, or that divide $F^*_q(0,1)$. This set is finite since $0$ is periodic. Enlarge the set $S$ if necessary so that the ring of $S$-integers $R_S$ is a PID. Write $y_k=r_k/s_k\in K$ with $r_k,s_k\in R_S$ and $(r_k,s_k)=(1)$. 
There are only finitely many elements $r_k$ and $s_k$, since $0$ is periodic.

Now write each $x_k=u_k/v_k\in K$ with $u_k,v_k\in R_S$ and $(u_k,v_k)=(1)$. Suppose that the prime $P$ is not in $S$ and that $m$ is the smallest nonnegative integer such that $P$ divides $F_m^*(u_0,v_0)$. Then $m$ is the smallest integer  such that $P$ divides $F_m(u_0,v_0)$, since $F_m(x,y)$ is (as we saw above) the product of   $F_{r}^*(x,y)$
to various powers, over $r\leq m$.

 Proceeding as in the proof of
Proposition 1, for $P \not\in S$ we see that 
	$$
		\phi^{(n)}(x_0) = x_n = \phi^{(n-m)}(x_m) \equiv \phi^{(n-m)}(0) = y_{n-m} \pmod P.
	$$
Hence $P$ divides $F_n(u_0,v_0)$
if and only if $P$ divides the numerator of $y_{n-m}$,
which holds if and only if $q$ divides $n-m$ as $P\not\in S$.
So if $P$ divides $F_n^*(u_0,v_0)$, we must have that $q$ divides $n-m$, and $n\geq m$.
If $n>m$, let $X=F_{n-q}(u_0,v_0), Y=G_{n-q}(u_0,v_0)$.  As $q$ divides
$n-m$ we know that $q$ divides
$(n-q)-m$, so that $P$ divides $F_{n-q}(u_0,v_0)=X$, and hence not $Y$ (else
$P$ divides Resultant$(F,G)$ which implies that $P\in S$, a contradiction). Now
$$F_n^*(u_0,v_0)=F_q^*(X,Y)\equiv F_q^*(0,Y) \equiv Y^\ell F^*_q(0,1) \pmod P,$$ where $\ell=\deg
F_q^*$.  But $P\nmid Y^\ell F^*_q(0,1)$ as $P\not \in S$, and so $P$ does not
divide $F_n^*(u_0,v_0)$.
Hence we have proved that if $P\not\in S$, then there is at most one value
of $n$ for which $P$ divides $F_n^*(x,y)$.

If $F_q(X,Y)$ has at least four non-proportional linear factors, then
$F_q(X,Y)/X^b$ has at least three non-proportional linear factors. Now
$F_n^*(x,y)=F_q(X,Y)/X^b$ where $X= F_{n-q}(x,y), Y=G_{n-q}(x,y)$, and the result
follows from the Thue/Mahler theorem, as in Proposition 1. (Note that the fact
that $\phi(t)\not\in \Cal T$ is part of the hypothesis.)

If $F_q(X,Y)$ has no more than three non-proportional linear factors,
select $k\geq 2$ minimal such that $F_{kq}(x,y)$ has at least
$2k+2$ non-proportional linear factors.  We know that such a $k$ exists
since Corollary 1 implies that  $F_{8q}(x,y)$ has at least
$d^{8q-4}+2\geq 2^{4}+2 = 2\cdot 8 +2$ non-proportional linear factors.
Then
$F_{kq}(X,Y)/F_{(k-1)q}(X,Y)^b$ has at least three non-proportional linear factors.
Now we see that  $F_n^*(x,y)=F_{kq}(X,Y)/F_{(k-1)q}(X,Y)^b$
where $X= F_{n-kq}(x,y), Y=G_{n-kq}(x,y)$, and the result
follows from the Thue/Mahler theorem, as in Proposition 1. \qed
\enddemo

\remark{Remark 4}
	Propositions~1 and~$1'$ imply the main results of [11], although no Diophantine approximation was necessary in the cases presented there. 
	\endremark

Proposition~$1'$ gives a unified means for finding prime numbers in certain residue classes.

\proclaim{Application} Let $q^n$ be an odd prime power. There exist infinitely many primes of the form $q^n k+1$.
\endproclaim

\demo{Proof}
	Consider the polynomial $\phi(t) = (t-1)^q + 1$. Then $\phi(0) =  0$ so that $0$ is a fixed point, but $\phi \not\in \Cal{T}$. Let $x_0$ be any integer larger than~1. Clearly $x_m \to \infty$ as $m\to \infty$, and so Proposition~$1'$ implies, after a small shift in notation, that $x_{n+m}$ has an odd primitive prime factor $p_m$ for all sufficiently large $m$. By the definition of $\phi$, we see that
$x_{m+r} - 1 = (x_{m+r-1} - 1)^q$ for every $r\geq 1$, and so by induction, 
	$$
		 (x_m - 1)^{q^n} = x_{m+n} -1 \equiv -1 \pmod{p_m}.
	$$
That is, $x_{m} - 1$ has order dividing $2q^n$ in the group $(\Bbb Z / p_m \Bbb Z)^{\times}$. 
If the order of $x_{m} - 1$ is $q^j$ for some $j \leq n$, then 
	$$
		-1 \equiv (x_{m} - 1)^{q^n} = \left((x_{m} - 1)^{q^j}\right)^{q^{n-j}} \equiv 1 \pmod{p_m},
	$$
a contradiction. On the other hand, if the order of $x_{m} - 1$ is $2q^j$, then
	$$
		(x_{m} - 1)^{q^j} \equiv -1 \pmod{p_m} \Rightarrow x_{m+j} = (x_{m}-1)^{q^j} +1 \equiv 0
			\pmod{p_m},
	$$
which contradicts the primitivity of $p_m$ unless $j = n$. Hence $x_{m} - 1$ has order exactly $2q^n$ in the group $(\Bbb Z / p_m \Bbb Z)^{\times}$, and consequently $2q^n $ divides $p_m - 1$. That is, $p_m \equiv 1 \pmod {q^n}$. Varying $m$, we produce infinitely many primes of this form.  \qed
\enddemo

\head 5. Dynamical systems with exceptional behavior at $\infty$, Part I \endhead

	In order to prove Theorem~1, we choose a point that falls into a cycle of length dividing $\Delta$ after exactly one step. The following lemma tells us this is always possible:
 
\proclaim{Lemma 2} Suppose $\phi(t) \in \Bbb C(t)$ has degree $d \geq 2$ and that $\Delta \geq 1$ is an integer. 
There exists a point $\alpha \in \Bbb P^1(\Bbb C) = \Bbb C \cup \{\infty\}$ such that $\phi^{(\Delta)}(\phi(\alpha)) = \phi(\alpha)$, but $\phi^{(\Delta)}(\alpha) \not= \alpha$.
\endproclaim

\demo{Proof} Suppose not. If $\beta$ is a fixed point of $\phi^{(\Delta)}$, and if
$\phi(\gamma)=\beta$ then $\gamma$ is also a fixed point of $\phi^{(\Delta)}$, or else we may take
$\alpha=\gamma$.  But then  
$\gamma=\phi^{(\Delta-1)}(\phi(\gamma))=\phi^{(\Delta-1)}(\beta)$ is unique,
and so $\gamma$ is totally ramified for $\phi$. By symmetry, $\beta$ is totally ramified too. 

In particular this implies that  $\phi'(\beta)=0$
and so  $(\phi^{(\Delta)}(\beta))'=\prod_{j=0}^{\Delta-1} \phi'(\phi^{(j)}(\beta))=0$. Hence $\beta$ has multiplicity 1 as a root of of $\phi^{(\Delta)}(x)-x$ by Lemma~A.1. Therefore 
$d^{\Delta} + 1$, the number of fixed points of $\phi^{(\Delta)}(x)-x$, by Lemma~A.2, 
equals the number of such $\beta$, which is
no more than the number of  totally ramified points of $\phi$. This is at most
$2$,  by the Riemann/Hurwitz formula (for any map $\phi$) and so we have established a   
contradiction. \qed
\enddemo

Define $\Cal E$ to be the class of rational functions $\phi(t) \in \Bbb C(t)$ of degree $d \geq2$ satisfying one of the following: 
	\roster
		\item"(i)" $ \displaystyle \phi(t) = t + \frac{1}{g(t)}$ for some polynomial $g(t)$ with $\deg (g) = d-1$;
		\item"(ii)" $\phi=\sigma^{-1}\circ \psi \circ \sigma$ for
some linear transformation $\sigma(t) = \lambda t+\beta$ with $\lambda\ne 0$, where $\displaystyle \psi(t)=\frac{t^d}{t^{d-1}+1}$; or
	\item"(iii)"  $\phi=\sigma^{-1}\circ \psi \circ \sigma$ with
 $\sigma$ as in (ii) and $\displaystyle \psi(t)=\frac{t^2}{2t+1}$.

	\endroster 

\noindent Observe that in the classes $\Cal E$(ii) and $\Cal E$(iii), the map $\sigma$ has no denominator, which means that $\infty$ is fixed under the transformation $\sigma$. The point at infinity plays an important role in the appearance of prime divisors in the numerators of our dynamical sequences.

\proclaim{Lemma 3} Let $\phi(t) \in \Bbb C(t)$ be a rational function of degree $d \geq 2$ and $\Delta$ a positive integer. If $\Delta \geq 2$, or if
$\Delta =1$ and $\phi(t) \not\in \Cal E$, then there exists 
$\alpha \in \Bbb P^1(\Bbb C)$ such that
$\phi^{(\Delta)}(\phi(\alpha))= \phi(\alpha),\ \phi^{(\Delta)}(\alpha) \not= \alpha$ and $\phi(\alpha) \not= \infty$.
\endproclaim

\demo{Proof} %We make a fractional linear change of coordinates to ensure that $\infty$ is not
%a root of $\phi^{(\Delta)}(x)$ or $\phi^{(\Delta+1)}(x)-\phi(x)$, and to send
%$\infty$ to $0$. Hence the lemma asks us to show that there exists   $\alpha \in  \Bbb C$ such that
%$\phi^{(\Delta)}(\phi(\alpha))= \phi(\alpha),\ \phi^{(\Delta)}(\alpha) \not= \alpha$ and $\phi(\alpha) \not= 0$.
We proceed as in the proof of Lemma 2:\  Suppose that the result is false.
If  $\beta\ne \infty$ is a fixed point of $\phi^{(\Delta)}$, then $\gamma = \phi^{-1}(\beta)$ is totally ramified for $\phi$, and it is a fixed point for $\phi^{(\Delta)}$. 
Moreover $\gamma$ is a fixed point of multiplicity one for $\phi^{(\Delta)}$ by Lemma~A.1. So if there are $r$ finite fixed points of $\phi$,
and $R$ finite fixed points of $\phi^{(\Delta)}$, then $r\leq R\leq 2$ by the 
Riemann/Hurwitz  formula (using the fact that the fixed points of $\phi$ are 
a subset of the fixed point of $\phi^{(\Delta)}$). If $\infty$ has multiplicity 
$k$ as a fixed point of $\phi$, then Lemma~A.2 gives  
	$$k+r=d+1.$$
We deduce that $k=d+1-r\geq 2+1-2=1$; that is,
$\infty$ must be a fixed point of $\phi$, and hence of $\phi^{(\Delta)}$.
% , which  $\phi^{-1}(\beta)\ne \infty$.
So either the finite fixed points of $\phi^{(\Delta)}$ are also fixed points of 
$\phi$, or there are two finite
fixed points $\beta, \gamma$ of $\phi^{(\Delta)}$ such that  
$\phi(\beta)=\gamma$ and  $\phi(\gamma)=\beta$. 
In either case we know that the finite fixed points of $\phi^{(\Delta)}$ are totally ramified points for $\phi$.

Now suppose that we have two finite fixed points $\beta, \gamma$ of $\phi^{(\Delta)}$
and change coordinates so that  $\beta \mapsto 0, \ \gamma \mapsto \infty$, and 
$\infty \mapsto 1$ to obtain a new function $\psi$. (When we say ``change coordinates so that $a \mapsto b$,'' we mean ``replace $\phi$ with $\psi = \sigma^{-1} \circ \phi \circ \sigma$,'' where $\sigma$ is a fractional linear transformation such that $\sigma^{-1}(a) = b$. Then $a$ is a fixed point of $\phi$ if and only if $b$ is a fixed point of $\psi$.) Now $\psi$
is totally ramified at $0$ and $\infty$, with pre-images $0$ and $\infty$,
and $\psi(1)=1$, so that $\psi(t)=t^{\pm d}$. As $\psi'(1) \not= 1$, we see that
$1$ is a fixed point of $\psi$ of multiplicity~1. Counting fixed points of $\psi$ using Lemma~A.2 and the
fact that $r$ is the number of fixed points distinct from 1 gives
$3\leq 1+r=d+1\leq 3$; that is, $r=d=2=R$. 
Now $\psi^{(\Delta)}(t)=t^{ 2^\Delta}$ (as $0$ and $\infty$ must be fixed), which has $2^{\Delta} + 1$ distinct fixed points, so that
$2^\Delta + 1 = 1+R=3$, and so $\Delta=1$. Hence the only possibility is
$\psi(t)=t^2$, and we obtain $\Cal E$(iii) as $\sigma^{-1}\circ \psi \circ \sigma$
with $\sigma(t)= t / (t+1)$. Note that the coordinate change at the beginning of this paragraph moved $\infty \mapsto 1$ and then this last coordinate change sent $1 \mapsto \infty$, so that $\infty$ was not moved by their composition. This corresponds to a change of coordinates of the form $\sigma(t) = \lambda t + \beta$ for some $\lambda \not= 0$.

Henceforth we may assume $R\leq 1$, so that the multiplicity of $\infty$ as a fixed point of $\phi$ satisfies $k=d+1-r\geq d+1-R\geq 2$. By Lemma~A.3(1), we know $\infty$ has multiplicity $k$ as a fixed point of $\phi^{(\Delta)}$, and so $d^\Delta+1$, the number of fixed points of $\phi^{(\Delta)}$, equals $k+R=d+1-r+R\leq d+2$. Hence $d(d^{\Delta-1}-1)\leq 1$, which implies that  $\Delta=1$.

If $R=0$, then $\infty$ is the only fixed point of $\phi$, which means $\phi(t) - t = 1 / g(t)$ for some polynomial $g(t)$ of degree $d-1$, from which we obtain  $\Cal E$(i).
 
If $R=1$, so that $\beta$ is the only finite fixed point of  $\phi$, replace $\phi$ with $\phi(t+\beta) - \beta$.  Now $\phi(t)=t^d/g(t)$, and the numerator of
$\phi(t)-t=t(t^{d-1}-g(t))/g(t)$ has only one root, so that $g(t)=t^{d-1}+ c$
for some constant $c \not= 0$. Taking $\psi(t)=\lambda^{-1} \phi(\lambda t)$ with $\lambda^{d-1} = c$ gives
$\Cal E$(ii). \qed

\enddemo

\head 6. Differences in the terms of dynamical sequences, Part I \endhead

	Theorem~1 is closely related to the following result, Theorem~$1'$, which gives primitive congruences in projective space. Their proofs are almost identical except that we need to be cautious about primes dividing the denominator in order to deduce Theorem~1. Theorem~$1'$ is perhaps more aesthetically appealing than Theorem~1 due to the fact that its conclusion holds for {\sl every} rational function.
	
\proclaim{Theorem $1'$} Suppose that $\phi(t)\in \Bbb Q(t)$ has degree $d\geq 2$,
and a positive  integer $\Delta$ is given. Let $x_0 \in \Bbb Q$ and define $x_{n+1} = \phi(x_n)$ for each $n \geq 0$.
If the sequence of rationals $\{ x_n\}_{n\geq 0}$ is not eventually periodic, then for all sufficiently large $n$ there exists a prime $p_n$ such that $x_{n+\Delta} \equiv x_n$ in $\Bbb P^1(\Bbb F_{p_n})$, but $x_{m+\Delta} \not\equiv x_m$ in $\Bbb P^1(\Bbb F_{p_n})$ for any $m < n$. 
\endproclaim

\demo{Proof of Theorem~$1'$, and then of Theorem~1}
	Choose a point $\alpha \in \Bbb P^1(\overline{\Bbb Q})$ such that $\phi^{(\Delta)}(\phi(\alpha)) = \phi(\alpha)$, but $\phi^{(\Delta)}(\alpha) \not= \alpha$, by Lemma~2. If   $\Delta >1$ or if $\phi(t) \not\in \Cal E$, then 
we  also insist that $\phi(\alpha) \not= \infty$, by Lemma~3.
Note that $\alpha$ is preperiodic, but not periodic for $\phi$. 
	
If $\alpha \in \Bbb C$ then define $\psi(t) := \phi(t+\alpha) - \alpha$,
else if $\alpha = \infty$  set $\psi(t) := 1 / \phi(1/t)$.
Note that $0$  is preperiodic, but not periodic for $\psi$ (and so 
$\psi(t) \not\in \Cal T$). Note that $\psi(t)\in K(t)$
for some finite Galois extension $K / \Bbb Q(\alpha)$. 

If $\alpha \in \Bbb C$ then let $y_0=x_0-\alpha$ and $y_{n+1}=\psi(y_n)$
for each $n\geq 0$;
one can easily verify that  $y_n=x_n-\alpha$ for all $n\geq 0$.
If $\alpha = \infty$ then let $y_0 = 1 / x_0$ and $y_{n+1} = \psi(y_n)$ for each $n\geq 0$;
now,  $y_n=1/x_n$ for all $n\geq 0$.

We apply Proposition~1 to the sequence $\{ y_n\}_{n\geq 0}$, so proving that the
numerator of $x_n-\alpha$ if $\alpha \in \Bbb C$, and the denominator of $x_n$
if $\alpha = \infty$, has a super-primitive prime factor $P_n$ (in $K$) for each sufficiently large $n$. By taking $n$
% a little 
larger if necessary, we may assume that
	
\roster
	\item"(i)"  No $P_n$ divides Resultant$(F,G)$ 
	
		\noindent (which guarantees that $\phi(t)$ induces a well-defined map of degree~$d$ on 
		$\Bbb P^1(\Bbb F_q)$ by canonically reducing each of its coefficients modulo $P_n$, 
		where $\Bbb F_q \cong R/P_nR$ is the finite field with $q$ elements); 

	\item"(ii)" No $P_n$ divides the denominator of $\phi(\alpha)$ if $\phi(\alpha)\ne\infty$; and
	
	\item"(iii)" $\phi^{(\Delta)}(\alpha) \not\equiv \alpha \pmod{P_n}$ for any $P_n$  
	
	\noindent (where the congruence is taken in $\Bbb P^1(R/P_nR)$, so that 
if $\alpha=\infty$ then condition (iii) means that $P_n$ does not divide the denominator of
$\phi^{(\Delta)}(\infty)$).
	
\endroster

We exclude only finitely many prime ideals in this way since $F(x,y)$ and $G(x,y)$ have no common linear factor over $\overline{\Bbb Q}$, and since $\alpha$ is not periodic.  

	The definition of $\alpha$ and the fact that $P_n$ divides the numerator of $x_n - \alpha$ yields
	$$
		x_{n+ 1+ \Delta} = \phi^{(\Delta + 1)}(x_n) \equiv \phi^{(\Delta + 1)}(\alpha)
		= \phi(\alpha) \equiv \phi(x_n) = x_{n+1} \pmod {P_n} 
	$$	
(i.e., in $\Bbb P^1(R/P_nR)$). If $p_n$ is the rational prime divisible by $P_n$  then 
$x_{n+1+\Delta} \equiv x_{n+1} \pmod {p_n}$, since $x_{n+1+\Delta}-x_{n+1}$ is rational.
Note that if  $\Delta >1$ or if $\phi(t) \not\in \Cal E$, then
$\phi(\alpha) \not= \infty$, and so, by condition (ii), $p_n$ does not divide the denominator of $x_{n+1}$ or $x_{n+ 1+ \Delta}$.

	We claim that $p_n$ is a primitive prime factor of $x_{n+1+\Delta} - x_{n+1}$. Indeed, suppose that $p_n$ is a factor of $x_{m+\Delta} - x_m$ for some $m < n+1$. Then
	$$
		\aligned
		\phi^{(\Delta)}(\alpha) \equiv \phi^{(\Delta)}(x_n) &= \phi^{(\Delta)}(\phi^{(n-m)}(x_m) ) \\
			&= \phi^{(n-m)}(x_{m+\Delta}) \equiv \phi^{(n-m)}(x_m) = x_n \equiv \alpha \pmod {P_n},
		\endaligned
	$$
contradicting the assumption in~(iii). We conclude that $p_n$ is a primitive prime factor of $x_{n+1+\Delta} - x_{n+1}$ and, changing variable to $N = n+1$, we deduce that there exists a primitive prime factor $p_n$ of $x_{N+\Delta} - x_N$ for all sufficiently large $N$.

This completes the proof of Theorem~$1'$. It also finishes the proof of Theorem~1 when $\Delta > 1$ or when $\Delta = 1$ and $\phi(t) \not\in \Cal E$, so it remains to treat the case $\Delta = 1$ and $\phi(t) \in \Cal E$. Let $x_n = u_n / v_n$ for coprime integers $u_n, v_n$. By Theorem~$1'$, either the numerator of $x_{n+1} - x_n$ or $v_n$, the denominator of $x_n$, has a primitive prime factor $p_n$.
	
Now if $\phi(t) = t + 1 / g(t)$ for some polynomial $g(t)$ of degree $d-1$, then
$$ 
x_{n+1} - x_n = \phi(x_n) - x_n = \frac{1}{g(x_n)} = \frac{v_n^{d-1}}{v_n^{d-1}g(u_n/v_n)}.
$$
Evidently the numerator is divisible by $p_n$, being a power of $v_n$, 
except perhaps if $p_n$ divides the leading coefficient of
$g$ (which can only occur for finitely many $n$). This completes the proof of Theorem~1 for such functions $\phi(t)$. 

	Similarly if $\displaystyle \lambda\phi(t)+\beta=   \frac{(\lambda t+\beta)^d}{(\lambda t+\beta)^{d-1}+1} $, so that $\phi \in \Cal E$(ii), then
$$
\lambda(x_{n+1} - x_n) = \lambda\phi(x_n) - \lambda x_n 
=  - \frac{(\lambda x_n+\beta)}{(\lambda x_n+\beta)^{d-1}+1}  =
- \frac{(\lambda u_n+\beta v_n)v_n^{d-2}}{(\lambda u_n+\beta v_n)^{d-1}+v_n^{d-1}},
$$
which is divisible by $p_n$ when $d > 2$ since $v_n$ is in the numerator, except perhaps if $p_n$ divides
the numerator of $\lambda$.  (Note that $\beta$ and $\lambda$ need not be rational, but the conclusion follows anyway upon consideration of prime ideal divisors.)

	For $d=2$ and  $\phi \in \Cal E$(ii), we study the function $\displaystyle \psi(t)=\frac{t^2}{t+1}$, proving that
$\displaystyle \psi^{(r+1)}(t)-\psi^{(r)}(t)= - \frac{t^{2^r}}{g_r(t)}$ where $g_r(t)$ is a monic polynomial
in $\Bbb Z[t]$ of degree $2^r$ by induction. For $r=0$ we have this with $g_0(t)=t+1$
by definition. For $r\geq 1$ we have, using the induction hypothesis,
$$
\psi^{(r+1)}(t)-\psi^{(r)}(t) = \psi^{(r)}(\psi(t))-\psi^{(r-1)}(\psi(t)) =
- \frac{\psi(t)^{2^{r-1}}}{g_{r-1}(\psi(t))} = 
- \frac{t^{2^r}}{ (t+1)^{2^{r-1}}g_{r-1}(\frac{t^2}{t+1})} = - \frac{t^{2^r}}{g_r(t)} .
$$
Hence the prime divisors of the numerator of $x_{r+1}-x_r$
are always the same: namely the prime divisors of the numerator of $x_0$.
Similarly,
for $\phi(t)$ obtained through the linear transformation $t\mapsto \lambda t+\beta$,
the prime factors of the numerator of $x_{r+1}-x_r$
are always the same, namely the prime divisors of the numerator of $\lambda x_0+\beta$ and
 the prime divisors of the denominator of $\lambda$.
 
Finally, we  show that the rational functions  $\Cal E$(iii) are exceptional:
For
$\psi(t)=\frac{t^2}{2t+1}$ we observe that 
$1+\frac{1}{\psi(t)}=\left(1+\frac 1t\right)^2$, and so
$1+\frac{1}{\psi^{(r)}(t)}=\left(1+\frac 1t\right)^{2^r}$
by an appropriate induction hypothesis. Hence
$$
\psi^{(r)}(t) = \frac 1 { \left(1+\frac 1t\right)^{2^r} - 1} \ \ \text{and} \ \ 
\psi^{(r+1)}(t)-\psi^{(r)}(t) = -\frac { \left(1+\frac 1t\right)^{2^r}} { \left(1+\frac 1t\right)^{2^{r+1}} - 1} =  \frac { \left( t(t+1) \right)^{2^r}} { t^{2^{r+1}} - \left(t+1\right)^{2^{r+1}}} .
$$
Therefore the prime divisors of the numerator of $x_{r+1}-x_r$
are always the same, namely the prime divisors of the numerator of $x_0(x_0+1)$. Similarly,
for $\phi(t)$ obtained through the linear transformation $t\mapsto \lambda t+\beta$,
the prime factors of the numerator of $x_{r+1}-x_r$
are always the same, namely the prime divisors of the numerator of $(\lambda x_0+\beta)(\lambda x_0+\beta+1)$ and
those of the denominator of $\lambda$. \qed
\enddemo

\remark{Remark 5} It is desirable to remove the Thue/Mahler Theorem from the above proof,
because, even though it is effective, the constants that come out are so large as to be of little
practical use. Moreover, the constants should grow with the field of definition of $\alpha$ (as chosen in the proof), and thus one should not expect any strong uniformity in $\Delta$ to come from this argument. So, do we really need the full power of the Thue/Mahler Theorem
in this proof?  In fact it may be the case that our proof can be modified 
to show that the exceptional $u_n$ must divide a particular non-zero
integer (rather than $u_n$ only having prime factors from a particular finite set). If we examine the proof above then we see that this idea works fine for the primes $P_n$ of types (ii)
and (iii). It is the primes that divide that resultant (i.e., those of type (i)) that require careful consideration
to determine whether their effect can be understood in this way. 
\endremark

\remark{Remark 6}
In the introduction we gave the example  $x_0=1$ with
$x_{n+1}=x_n^2/(2x_n+1)$ so that $x_n = \frac 1{F_n-2}$, and
$x_{n+1}-x_n=-\frac{2^{2^n}}{F_{n+1}-2}$.
Another amusing example is given by 
%``Euclid numbers'' 
Sylvester's sequence $E_0=2$ and $E_{n+1}=E_n^2-E_n+1$. (The terms are $2, 3, 7, 43, \ldots$, which can occur
in a version of Euclid's proof of the infinitude of primes, based on the fact that
$E_n=E_{n-1}E_{n-2}\cdots E_0+1$).  Now let $x_0=1$ with
$x_{n+1}=x_n^2/(x_n+1)$ so that $x_n = \frac 1{E_n-1}$, and
$x_{n+1}-x_n=-\frac{1}{E_n}$, so that there are never prime divisors of
the numerator of $x_{n+1}-x_n$.
\endremark

\remark{Remark 7}
	One could instead use Proposition~$1'$ to prove Theorem~1. In this case, we would have to select a point $\alpha \in \Bbb P^1(\overline{\Bbb Q}) \smallsetminus \{\infty\}$ that is periodic with period dividing $\Delta$, but that is not totally ramified. This allows us to change coordinates to obtain a new rational function that is not in $\Cal T$. One can choose a fixed point $\alpha$ with this property precisely when $\phi \not\in \Cal E$, and then the proof proceeds essentially as above.
\endremark

\remark{Remark 8} If $p$ is a prime dividing the numerator of $x_{n+\Delta}-x_{n}$ but not Resultant$(F,G)$, then
$$x_{n+1+\Delta}=\phi(x_{n+\Delta})\equiv \phi(x_{n})=x_{n+1} \pmod {p},$$
and so $p$ divides the numerator of $x_{n+1+\Delta}-x_{n+1}$. Iterating we find that
$p$ divides the numerator of $x_{m+\Delta}-x_{m}$ for all $m>n$.

We can also ask to understand the power of $p$ appearing as a factor in each subsequent term. From the Taylor expansion
$\phi(t+h)=\phi(t)+h\phi'(t)+\frac {h^2}{2!} \phi''(t)+\cdots$
with $h=\phi^{(\Delta)}(t)-t$, we deduce that
$\phi^{(\Delta+1)}(t)-\phi(t) - (\phi^{(\Delta)}(t)-t)\phi'(t)$ is divisible by
$(\phi^{(\Delta)}(t)-t)^2$. Taking $t=x_n$ we deduce that
$x_{n+\Delta}-x_{n}$ divides $x_{n+1+\Delta}-x_{n+1}$, up to a bounded quantity,
so we recover the result of the previous paragraph. But we can go much further
assuming that  the numerators
of $\phi^{(\Delta)}(t)-t$ and $\phi'(t)$ have no common factor. If so, then
the gcd of the numerators of $\phi^{(\Delta)}(x_n)-x_n$ and $\phi'(x_n)$
divides the resultant of the  polynomials in the two numerators, which is
non-zero. Hence for all but finitely many primes $p_n$, if $p_n^e\| x_{n+\Delta}-x_{n}$
where $p_n\nmid x_m$ for all $m<n$,  then $p_n^e\| x_{N+\Delta}-x_{N}$
for all $N\geq n$.
\endremark

\head 7. Baker's Theorem and primitive prime factors  \endhead

\proclaim{Baker's Theorem}([1, Thm.~3]) A rational map of degree $d \geq 2$ defined over $\Bbb{C}$ has a periodic point in $\Bbb{P}^1(\Bbb{C})$ of {\sl exact} period $\Delta \geq 2$ except perhaps when $\Delta=2, d=2,3$ or $4$, or when $\Delta=3, d=2$. There exist exceptional maps in each of these four cases. 
\endproclaim

The exceptional maps in Baker's theorem are also exceptions to our Theorem~$2$, as is shown by Lemma~4 below. The exceptions were classified up to conjugation by a linear fractional transformation by Kisaka [7]; see Appendix~B for the classification. As in the introduction, let us write $\Cal B_{\Delta, d}$ for the set of rational functions of degree~$d$ with no point of exact period $\Delta$; so $\Cal B_{2,2} \cup \Cal B_{2,3} \cup \Cal B_{2,4} \cup \Cal B_{3,2}$ is the set of exceptions in Baker's Theorem.

\proclaim{Lemma 4}
	Suppose $\phi(t) \in \Cal B_{\Delta,d} \cap \Bbb Q(t)$ is a rational map of degree~$d$ with no periodic point of period $\Delta$. There exists a finite set of primes $S$ such that for any $u/v \in \Bbb Q$, either $\phi(u/v)= \infty$, $\phi^{(\Delta)}(u/v) = \infty$, or else every $p \not\in S$ that divides the numerator of  $\phi^{(\Delta)}(u/v) - u/v$ is also a factor of the numerator of $\phi(u/v) - u/v$. 	
\endproclaim

\demo{Proof}
	Write
		$$
			\aligned
				\phi(x/y) - x/y &= \frac{A_1(x,y)}{B_1(x,y)} \\
				\phi^{(\Delta)}(x/y) -x/y &= \frac{A_\Delta(x,y)}{B_\Delta(x,y)},
			\endaligned
		$$
where $A_i(x,y), B_i(x,y) \in \Bbb Z[x,y]$ are homogeneous polynomials such that $A_1$ and $B_1$ (resp. $A_\Delta$ and $B_\Delta$) share no common linear factor over $\overline{\Bbb Q}$ and no common factor in their content. (Recall that the content of a polynomial with integer coefficients is the greatest common divisor of its coefficients.) As $\phi(t)$ has no point of exact period $\Delta$, every solution in $\Bbb{P}^1(\overline{\Bbb{Q}})$ to $\phi^{(\Delta)}(\alpha) = \alpha$ must also be a solution to $\phi(\alpha) = \alpha$. (Here we are using the fact that $\Delta = 2$ or $3$ is prime.) In particular, any non-constant factor of $A_\Delta(x,y)$ that is irreducible over $\Bbb Z$ is also a factor of $A_1(x,y)$. 

	Define $S$ to be the set of primes dividing Resultant$(A_1, B_1)$ together with those primes dividing the content of $A_\Delta$. We may assume that $\phi(u/v)$ and $\phi^{(\Delta)}(u/v)$ are not equal to infinity, and also that $u,v$ are coprime integers. Suppose $p \not\in S$ is a prime factor of the numerator of $\phi^{(\Delta)}(u/v) - u/v$. Then $p \mid A_\Delta(u,v)$, and consequently there exists an irreducible factor (over $\Bbb Z$) of $A_\Delta$, say $Q(x,y)$, such that $p \mid Q(u,v)$. By the last paragraph, we know $Q(x,y)$ divides $A_1(x,y)$, and hence $p \mid A_1(u,v)$. Now $p \nmid B_1(u,v)$ since otherwise $p \mid $ Resultant$(A_1, B_1)$. We conclude that $p$ divides the numerator of $\phi(u/v) - u/v$. \qed
\enddemo	

\proclaim{Corollary 2} Suppose $(\Delta,d)$ is one of the exceptional pairs in Baker's Theorem, and let $\phi(t) \in \Cal B_{\Delta,d} \cap \Bbb Q(t)$. There is a finite set of primes $S$ with the following property. If we define $x_0 \in \Bbb Q$ and $x_{n+1} =\phi(x_n)$, and if the sequence $\{x_n\}_{n \geq 0}$ is not eventually periodic, then for all $n$ sufficiently large, any prime $p \not\in S$ that divides the numerator of $x_{n + \Delta} - x_n$ must also divide the numerator of $x_{n+1} - x_n$. 

\endproclaim

\head 8. Dynamical systems with exceptional behavior at infinity, Part II \endhead

	Recall that $p_{\Delta, n}$ is a doubly primitive prime factor of $x_{n+\Delta} - x_n$ if $p_{\Delta, n}$ divides the numerator of $x_{n+\Delta} - x_n$, and if $N\geq n$ and $D\geq \Delta$ whenever 
$p_{\Delta, n}$ divides the numerator of $x_{N+D} - x_N$. To produce a doubly primitive prime factor of the numerator of $x_{n+\Delta} - x_n$, we want an $\alpha \in \Bbb P^1(\overline{\Bbb Q})$ such that $\phi(\alpha)$ is not $\infty$, and $\phi(\alpha)$ has exact period $\Delta$, while $\alpha$ is not itself periodic. This will allow us to apply Proposition~1 inductively as in the proof of Theorems~1 and~$1'$.

\proclaim{Lemma 5} Suppose $\phi(t) \in \Bbb C(t)$ is a rational function of degree $d \geq 2$, and let $\Delta \geq 1$ be an integer. There exists a point $\alpha \in \Bbb P^1(\Bbb C) = \Bbb C \cup \{\infty\}$ such that $\phi(\alpha)$ has exact period $\Delta$, $\phi(\alpha) \not= \infty$, and $\phi^{(\Delta)}(\alpha) \not= \alpha$ unless 
	\roster
		\item"(1)" $\Delta = 1$ and $\phi(t) \in \Cal E$ (see \S5 for the definition of $\Cal E$); or
		
		\item"(2)" $\Delta = 2$ and $\phi(t) = \alpha + 1 / g(t-\alpha)$ for some $\alpha \in \Bbb C$ and 
			some quadratic polynomial $g(t) \in \Bbb C[t]$ such that $g(0) = 0$, but $g(t) \not= ct^2$ 
			for any complex number $c$; or
		
		\item"($2'$)" $\Delta = 2$ and $\phi = \sigma^{-1} \circ \psi \circ \sigma$ for some 
			$\sigma(t) = \lambda t + \beta$ with $\lambda \neq 0$, and $\psi(t) = 1 / (t^3 +2t^2 + 2t)$; or
			
		\item"(3)" $\Delta = 2$ and $\phi = \sigma^{-1} \circ \psi \circ \sigma$ for some 
			$\sigma(t) = (\alpha t + \beta) / (\gamma t + \delta) $ with $\alpha \delta - \beta \gamma \not= 0$, and
				 $\psi(t) = 1/t^2$; or
			
		\item"(4)" $\Delta = 2$ and $\phi \in \Cal B_{2, d}$ some $d = 2, 3, 4$ 
			(see \S7 for the definition of $\Cal B_{\Delta, d}$); or
		
		\item"(5)" $\Delta = 3$ and $\phi \in \Cal B_{3,2}$
	\endroster
\endproclaim

\demo{Proof}  Assume that no such $\alpha$ exists.
Baker's Theorem states that if $\phi$ has no point of exact period $\Delta$ then we
are in case  (4) or (5). So henceforth assume  that $\phi$ has a point $\gamma$ of exact period $\Delta$. 
Every point in the orbit of $\gamma$ must also have exact period $\Delta$, and hence there are at least $\Delta$ points of exact period $\Delta$.

As in the proofs of Lemmas 2 and 3, if $\beta\ne \infty$ has period $\Delta$, then 
$\gamma = \phi^{-1}(\beta)$ must be totally ramified for $\phi$ else there would be 
an $\alpha$ as desired with $\phi(\alpha)=\beta$.
There are no more than two elements
that are  totally ramified for (any rational map) $\phi$,
hence there can be no more than  three points of
exact period $\Delta$ (that is, $\infty$ and the two points with fully ramified pre-images). Hence
$\Delta\leq 3$.

The case $\Delta = 1$ is given by Lemma~3 in \S5: The exceptions are precisely those in $\Cal E$;
that is, case (1). 
 
If $\Delta=2$ or $3$ then all of the points of exact order $\Delta$ must be in a unique
orbit, else there would be at least $2\Delta\geq 4$ points of exact order $\Delta$, a contradiction.

So, if $\Delta=3$ then  there is exactly one orbit, and it is of the form $\infty \mapsto \beta_1 \mapsto \beta_2 \mapsto \infty$  with $\infty$ and $\beta_1$ being totally ramified. We will show this is impossible. 
% $\infty$ and two totally ramified points which, we will show, is impossible. 
We may conjugate by a linear fractional transformation in order to assume the totally ramified points are $1$ and $2$, and that the other fixed point is $0$. We suppose that $\phi(2)=1, \phi(1)=0,\ \phi(0)=2$. The only possible
rational function with these ramification conditions is
$$
	\phi(t)  = \frac{ 2^{d+1} (t-1)^d} { 2^{d+1} (t-1)^d-(t-2)^d},
$$
so that
$$
\phi^{(3)}(t) - t = \phi(\phi^{(2)}(t))-t = \frac{(1-t) 2^{d+1} (\phi^{(2)}(t)-1)^d  + t(\phi^{(2)}(t)-2)^d}{2^{d+1} (\phi^{(2)}(t)-1)^d-(\phi^{(2)}(t)-2)^d}.
$$
We choose homogeneous polynomials $F_2(x,y), G_2(x,y) \in \Bbb C[x,y]$ of degree~$d^2 = \deg \phi^{(2)}$ with no common linear factor so that $\phi^{(2)}(x/y) = F_2(x,y) / G_2(x,y)$ after clearing denominators. (One can use the polynomials defined in \S2, for example.) The numerator of $\phi^{(3)}(x/y)-x/y$ is therefore
$$
N_3(x,y)  =  (y-x) 2^{d+1} \left(F_2(x,y)-G_2(x,y)\right)^d  + x\left(F_2(x,y)-2 G_2(x,y)\right)^d.
$$
Dividing through by the common factors $x(x-y)$ we will apply the $abc$-theorem for polynomials (\S3).  Note that the number of non-proportional linear factors of $N_3(x,y)$ is
$3$ plus the number of distinct fixed points of $\phi(t)$, which is at most
$3+(d+1)=d+4$.  Hence the total number of distinct roots
in our $abc$-equation is at most $2\deg \phi^{(2)} + d+4$. The $abc$-theorem implies that
$(d\deg \phi^{(2)}-1) + 2 \leq 2\deg \phi^{(2)} + d+4$, hence
$d^2(d-2)\leq d+3$, and so $d=2$.  When $d=2$ we find that the numerator of $\phi^{(3)}(t) - t$ is 
$$
	t(t-1)(t-2)(7t^3-14t^2+8)N_1(t),
$$	
where $N_1(t)$ is the numerator of $\phi(t) - t$. It follows that we have a second
cycle of exact order three, consisting of the roots of $7t^3-14t^2+8$. But this contradicts our hypothesis that only one orbit of length~3 exists. 

We know that if $\Delta = 2$ then there is a single orbit of length~2. If
both points of order $2$ are totally ramified, we change coordinates so that 
the $2$-cycle consists of $0$ and $\infty$, and that $1$ is a fixed point. 
Then $\psi(t) = 1/t^d$, and hence
	$$
		\psi(t) - t = (1-t^{d+1})/{t^d}; \qquad \psi^{(2)}(t)-t = t (t^{d^2-1}-1).
	$$
Since there are no points of exact order 2 other than $0$ and $\infty$, all of the
$(d^2-1)^{\operatorname{th}}$ roots of unity (which satisfy $\psi^{(2)}(t)=t$) must also
be $(d+1)^{\operatorname{th}}$ roots of unity, so as to satisfy $\psi(t)=t$. Hence $d^2 - 1\leq d+1$, 
and thus $d=2$  and $\psi(t) = 1 / t^2$. It follows that any rational map of the form $\phi(t) = (\sigma^{-1} \circ \psi \circ \sigma)(t)$ with $\sigma(t) = (\alpha t+\beta) / (\gamma t + \delta)$ and $\alpha \delta - \beta \gamma \not= 0$ has a unique periodic orbit of length~2 consisting entirely of totally ramified points,
yielding case (3).

	Finally suppose that $\Delta = 2$ and there is a unique cycle of length~2 consisting of
$\infty$ --- which is totally ramified --- and one other point $\alpha \in \Bbb C$. We may assume $\alpha$ is not totally
ramified, else we are in the previous case. We have
$\phi(t) = \alpha + 1 / g(t-\alpha)$ where $\deg(g) = d$ and $g(0) = 0$, and $g(t) \not= ct^d$
for some complex number $c$. The remainder of the proof follows the strategy of Baker's Theorem. The map $\phi^{(2)}$ has $d^2+1$ fixed points counted with multiplicity (Lemma~A.2). Two of them are $\alpha$ and $\infty$, which each have multiplicity $1$ since the fixed point multiplier of $\phi^{(2)}$ at $\alpha$ and $\infty$ is zero 
(because $(\phi^{(2)})'(\alpha)=(\phi^{(2)})'(\infty)=\phi'(\infty)\phi'(\alpha)=0$ as $\infty$ is ramified). All of the remaining fixed points of $\phi^{(2)}$ must be fixed points of $\phi$ as we are assuming there are no other periodic orbits of length~2. By Lemma~A.3, each of the remaining fixed points $\beta_i \in \Bbb C$ falls into exactly one of the following categories:
	\roster
		\item"(i)" The fixed point multiplicity of $\beta_i$ for $\phi$ and $\phi^{(2)}$ is 1. Suppose there are $M$
			fixed points of this type.
			
		\item"(ii)" The fixed point multiplicity of $\beta_i$ for $\phi$ is $\ell_i > 1$, which implies the fixed point
			 multiplicity for $\phi^{(2)}$ is also $\ell_i$. 
			 
		\item"(iii)" The fixed point multiplicity of $\beta_i$ for $\phi$ is $1$, and the fixed point multiplicity of 
			$\beta_i$ for $\phi^{(2)}$ is $2k_i + 1$ for some positive integer $k_i$. Suppose there are $r$
			fixed points of this type.
	\endroster
As $\phi$ has exactly $d+1$ fixed points (always with multiplicity), we may add up the fixed points of these types to find
	$$
		d+1 = M + \sum \ell_i + r.
	$$
Applying the same reasoning to $\phi^{(2)}$ and noting that we must also count $0$ and $\infty$, we have
	$$
		d^2 + 1 = 2 + M + \sum \ell_i + \sum (2k_i + 1) = 2 + M +\sum \ell_i + r +  2\sum k_i.
	$$
Subtracting the first of these equations from the second gives $d^2 - d = 2 + 2\sum k_i$, or  $\sum k_i = \frac{1}{2}(d^2 - d - 2)$. On the other hand, Lemma~A.3 also tells us that each of the type (iii) fixed points $\beta_i$ attracts $k_i$ distinct critical points. Since a critical point can only be attracted to a single one of the $\beta_i$, we see that there are $\sum k_i$ distinct critical points attracted to the set of fixed points $\{\beta_i\}$. By the Riemann/Hurwitz formula, there are exactly $2d-2$ critical points (with multiplicity). Now $\infty$ is a critical point of order $d-1$ (as it is totally ramified), but it is also a periodic point, so it cannot be attracted to one of the $\beta_i$. Hence 
	$$
		\frac{1}{2}(d^2 - d - 2) = \sum k_i \leq d-1 \Longrightarrow d \leq 3.
	$$
	
	 If $d = 2$, then $\phi$ has three fixed points counted with multiplicity, and $\phi^{(2)}$ has $5$ (Lemma~A.2). Hence there are exactly two points of exact period~$2$, namely $\alpha$ and $\infty$. This yields case (2). 
	 
	 If $d = 3$, let us change coordinates so that $\alpha = 0$. Now $\phi(t) = 1 / g(t)$, where $g(t) = at^3 + bt^2 + ct$ and $a\not= 0$. Choosing $\delta \in \Bbb C$ such that $\delta^4 = a^{-1}$ and replacing $\phi(t)$ by $\delta^{-1} \phi(\delta t)$, we may even suppose that $a = 1$. A direction calculation shows the numerator of $\phi(t) - t$ is
	 $$
	 	t^4+b t^3+c t^2-1,
	 $$
while the numerator of $\phi^{(2)}(t) - t$ is
	$$
		t(t^4+b t^3+c t^2-1)\left(t^4+2 b t^3-(-b^2-c) t^2+b c t+1\right).
	$$
Therefore $\phi$ has a second periodic orbit of length~2 if the polynomials $t^4+b t^3+c t^2-1$ and $t^4+2 b t^3-(-b^2-c) t^2+b c t+1$ have no common root. The resultant of these two polynomials is $b^4 - 4b^2c+16$, which shows that they have a common root if and only if $b \not= 0$ and $c = (b^4+16) / 4b^2$. Let us now assume that $c = (b^4+16)/ 4b^2$, in which case the numerators of $\phi(t) - t$ and $\phi^{(2)}(t) - t$ become
	$$
		(2t+b)(2b^2 t^3 + b^3 t^2 +8 t-4b)  \text{ and } 
		t(2t+b)^3(2b^2 t^3 + b^3 t^2 +8 t-4b)(4b^4 t^2  +4 t b^5 +16 b^2).
	$$
The roots of the final factor $4b^4 t^2  +4 t b^5 +16 b^2$ will yield a new periodic orbit of length~2 unless this polynomial shares a root with the numerator of $\phi(t)  - t$. The resultant of these two polynomials is
	$$
		-4b^{12}(b^4 - 16).
	$$
The case $b = 0$ has been ruled out by an earlier part of the argument, so we must have $b = 2 i^m$ for some $i = 0, 1, 2, 3$. (Here $i^2 = -1$.) Thus $c = (b^4 + 16) / 4b^2 = 2 (-1)^m$, and hence
	$$
		\phi(t) = \frac{1}{t^3 + 2 i^m t^2 + 2 (-1)^m t}, \quad m = 0, 1, 2, 3.
	$$
Replacing $\phi$ with $i^{-m}\phi(i^m t)$ yields the map $\psi$ in case ($2'$). 
   \qed
\enddemo

Lemma 5 is used in the proof of Theorem 2 as in the following:

\proclaim{Lemma 6} Suppose that $\phi(t) \in \Bbb Q(t)$ is a rational function of degree $d \geq 2$, and 
that there exists a point $\alpha \in \Bbb P^1(\overline{\Bbb Q})$ such that $\phi(\alpha)$ 
has exact period $\Delta\geq 1$, $\phi(\alpha) \not= \infty$, and $\phi^{(\Delta)}(\alpha) \not= \alpha$.
Let $x_0 \in \Bbb Q$ and define $x_{n+1} = \phi(x_n)$ for each $n \geq 0$, and suppose
that  the sequence $\{ x_n\}_{n\geq 0}$ is not eventually periodic. Suppose that
$P$ is a prime ideal that does not divide {\rm Resultant}$(f,g)$ (where $\phi=f/g$), and that
$p$ is the rational prime divisible by $P$.
If $P$ divides the numerator of $x_n-\alpha$, but neither the denominator of $\alpha$ nor $\phi(\alpha)$, and neither the numerator of $\phi^{(\Delta)}(\alpha)- \alpha$ nor
$\phi^{(\ell)}(\phi(\alpha))-\phi(\alpha)$ for any $1\leq \ell<\Delta$, then the
prime $p$ divides the numerator of $x_{N+D}-x_N$ if and only if $N\geq n+1$ and $\Delta$ divides $D$.
\endproclaim

\demo{Proof}  We begin by noting that 
$$
x_{n+\Delta+1}=\phi^{(\Delta)}(\phi(x_n))\equiv \phi^{(\Delta)}(\phi(\alpha)) = \phi(\alpha) \equiv 
\phi(x_n)=x_{n+1} \pmod P,
$$
and so $p$ divides the numerator of $x_{n+1+\Delta}-x_{n+1}$. (We have used the fact that 
$P$ does not divide the denominator of $\phi(\alpha)$.)

If $N\equiv n+j \pmod \Delta$ for $1\leq j\leq \Delta$ with $N>n$ then
$x_N\equiv x_{n+j} \pmod p$. To see this, we proceed by induction on
$N\geq n+1+\Delta$ since
$$
x_N=\phi^{(N-({n+1+\Delta}))}(x_{n+1+\Delta}) \equiv 
\phi^{(N-({n+1+\Delta}))}(x_{n+1}) = x_{N-\Delta} \pmod p .
$$

We now prove that if $p$ divides the numerator of $x_{N+D}-x_N$, then $\Delta$ divides $D$. If not let
$D\equiv \ell \pmod \Delta$ where $1\leq \ell<\Delta$, and select $m$ a large integer
%$m > N+D$
such that 
$m\equiv n+1 \pmod \Delta$. Then, using the congruences of the previous paragraph,
$$
x_{n+1} \equiv x_m = \phi^{(m-N)}(x_N) \equiv \phi^{(m-N)}(x_{N+D}) = x_{m+D}\equiv 
x_{n+1+\ell}  \pmod p.
$$
Hence
$$
\phi(\alpha)\equiv \phi(x_n)=x_{n+1} 
\equiv x_{n+1+\ell}  = \phi^{(\ell+1)}(x_n) = \phi^{(\ell)}(\phi(\alpha))  \pmod P,
$$
which contradicts the hypothesis.

Finally suppose that $p$ divides the numerator of $x_{N+D}-x_N$ with $N\leq n$ and $\Delta$ divides $D$. Then
$$
x_n=\phi^{(n-N)}(x_N) \equiv \phi^{(n-N)}(x_{N+D}) = x_{n+D}\equiv 
x_{n+\Delta}  \pmod p ,
$$
using the congruence from two paragraphs above, and so 
$$
\alpha\equiv  x_n \equiv x_{n+\Delta} 
=   \phi^{(\Delta)}(x_n) \equiv \phi^{(\Delta)}(\alpha)  \pmod P 
$$
which contradicts the hypothesis.  \qed

\enddemo

\head 9. Differences in the terms of dynamical sequences, Part II \endhead

	Now we give the proof of Theorem~2. We recall the statement for the reader's convenience. Define $\Cal F_1$ to be those $\phi(t)\in \Bbb Q(t)$ of the form 
$\sigma^{-1}\circ \psi \circ \sigma$, for
some linear transformation $\sigma(t) = \lambda t+\beta$ with $\lambda\ne 0$, where
	$$
		\psi(t) = \frac{t^2}{t+1} \ \text{or} \ \frac{t^2}{2t+1} .
	$$
Define $\Cal F_2$ to be the union of $\Cal B_{2, d}$ for $d = 2, 3, 4$ along with all rational maps $\phi(t)$ of the form $\phi = \sigma^{-1} \circ \psi \circ \sigma$ for some $\sigma(t) = (\alpha t - \beta) / (\gamma t - \delta)$ with $\alpha \delta - \beta \gamma \not= 0$ and $\psi(t) = 1 / t^2$. Define $\Cal F_3 = \Cal B_{3,2}$. 

\proclaim{Theorem 2} Suppose that $\phi(t)\in \Bbb Q(t)$ has degree $d\geq 2$.
Let $x_0 \in \Bbb Q$ and define $x_{n+1} = \phi(x_n)$ for each $n \geq 0$, and suppose
that  the sequence $\{ x_n\}_{n\geq 0}$ is not eventually periodic. For any given 
$M\geq 1$, the numerator of $x_{n+\Delta}-x_n$ has a doubly primitive prime factor for 
all $n \geq 0$ and $M \geq  \Delta \geq 1$, except  for
those pairs with $\Delta = 1, 2$ or $3$ when $\phi \in \Cal F_1, \Cal F_2$ or $\Cal F_3$, respectively,
as well as for finitely many other exceptional pairs $(\Delta, n)$.
 \endproclaim

\demo{Proof of Theorem 2} We proceed by induction on $M$, the case $M = 1$ being a consequence of Theorem~1. Suppose now that the result holds for all $M \leq M'-1$, and let us prove it holds for $M = M' \geq 2$. By the induction hypothesis, we  find that the numerator of $x_{n+\Delta} - x_n$ has a doubly primitive prime factor $p_{\Delta,n}$ for $n \geq 0$ and $M' - 1 \geq \Delta \geq 1$ other than for finitely many pairs $(\Delta, n)$, excluding those with $\Delta = 1, 2$ or $3$ and $\phi \in \Cal F_1, F_2$ or $F_3$, respectively. Set $\Delta = M'$. 

	Let us suppose that $\phi$ does not belong to one of the corresponding exceptional classes of Lemma~5, in which case we can choose $\alpha \in \Bbb P^1(\overline{\Bbb Q})$  so that $\phi(\alpha)$ has exact period $\Delta$, $\phi(\alpha) \not= \infty$, and $\phi^{(\Delta)}(\alpha) \not= \alpha$. Proceeding as in the proof of Theorems~$1'$ and~1, we obtain a sequence of prime ideals $\{P_{\Delta,n}\}$ of the Galois closure of $\Bbb Q(\alpha)$ such that $P_{\Delta,n}$ is a primitive prime factor of the numerator of $x_n - \alpha$ for all sufficiently large $n$. Moreover, if $p_{\Delta,n}$ is the rational prime divisible by $P_{\Delta,n}$, then $p_{\Delta,n}$ is a doubly primitive prime factor of the numerator of $x_{n+1+\Delta} - x_{n+1}$ by Lemma 6, provided that $P_{\Delta,n}$ is not one of the finitely many
prime ideal divisors of the numerator of $\phi^{(\Delta)}(\alpha)- \alpha$ or  $\phi^{(\ell)}(\phi(\alpha))-\phi(\alpha)$ for some $1 \leq \ell \leq \Delta - 1$, or of the denominator of $\alpha$ or~$\phi(\alpha)$.
	
	 Now suppose that $\phi$ belongs to case (2) of Lemma~5, so that $\Delta = M' = 2$ and that $\phi(t) = \alpha + 1 / g(t-\alpha)$ for some $\alpha \in \overline{\Bbb Q}$ and some quadratic polynomial $g(t) = bt^2 + ct$ with $bc \not= 0$. Define $\psi(t) = \phi(t+\alpha) - \alpha = 1 / g(t)$, so that $0$ is periodic of period~2 for $\psi$ and $\psi(0) = \infty$. Note $\psi \not\in \Cal T$ (as in \S3) since $g(t) \not= ct^2$. Define $y_0 = x_0 - \alpha$ and $y_{n+1}  = \psi(y_n)$. Then $y_n = x_n - \alpha$ by induction. Let $K/\Bbb Q$ be a Galois extension containing $\Bbb Q(\alpha)$. Invoking Proposition~$1'$, we see there exists a prime ideal $P_{2, n}$ of $K$ that is a primitive prime factor of the numerator of $y_n = x_n - \alpha$ for all sufficiently large $n$. We exclude the finitely many prime ideals
$P_{2, n}$ that divide the numerator or denominator of $b$ or $c$.	
Note that, since $\psi(x/y)= y^2 /x (bx + cy)$ , $P_{2, n}$ is a primitive
prime factor of the denominator of $y_{n-1}$, and hence for $N\geq n-1$, we have that 
$P_{2, n}$ divides the numerator of $y_N$ if $N-n$ is even, and the denominator of $y_N$ if $N-n$ is 
odd. Moreover,  if $e$
is the exact power of $P_{2, n}$ dividing the denominator of $y_{n-1}$, then 
$2^le$ is the exact power of 
$P_{2, n}$ dividing the numerator of $y_{n+2l-2}$ and  the denominator of $y_{n+2l-1}$
for all $l\geq 1$.
Hence if $a>b\geq n-1$ then $P_{2, n}$ divides the numerator of $x_a-x_b=y_a-y_b$ if and only
if  $a\equiv b\equiv n \pmod 2$.
Now suppose that $P_{2, n}$ divides the numerator of $x_a-x_b=y_a-y_b$ with $b<n-1$. If $a\leq n-1$ we have 
$y_{n-1}=\psi^{(n-1-a)}(y_a)\equiv \psi^{(n-1-a)}(y_b)=y_{n-1-(a-b)} \pmod {P_{2, n}}$, so we may
assume, without loss of generality that $a\geq n-1$.  Therefore $P_{2, n}$ divides the numerator
or denominator of $y_a$, and so of $y_b$ (as $P_{2, n}$ divides the numerator of their difference),  which contradicts
primitivity. In summary, we have shown that the numerator of $x_{N + D} - x_N$ is divisible by $P_{2,n}$ if and only if $N \geq n$ and $D \geq 2$ is even, and since $x_{N+D} - x_N$ is rational, the same statement is true when $P_{2,n}$ is replaced by the rational prime $p_{2,n}$ dividing $P_{2,n}$. This completes the proof for the maps from case (2) of Lemma~5.

	The strategy of the previous paragraph carries over {\it mutatis mutandis} to case ($2'$) of Lemma~5.

	Now suppose that $\phi$ belongs to case (3) of Lemma~5, so that  $\Delta = M' = 2$ and there exists a fractional linear tranformation $\sigma(t) = (\alpha t + \beta)/(\gamma t + \delta)$ such that $\phi = \sigma^{-1} \circ \psi \circ \sigma$ with $\psi(t) = 1/t^2$. Set $K = \Bbb Q(\alpha, \beta, \gamma, \delta)$, let $R$ be the ring of integers of $K$, and let $S$ be a finite set of prime ideals of $R$ such that the ring of $S$-integers $R_S$ is a principal ideal domain, and such that $\alpha, \beta, \gamma, \delta \in R_S$. (In particular, $S$ contains all prime ideals dividing the denominators of $\alpha, \beta, \gamma, \delta$.) Let $y_0=\sigma(x_0)$ and $y_{n+1}=\psi(y_n)$ for all $n\geq 0$. For each $n$ choose $u_n, v_n \in R_S$ such that $y_n = u_n/v_n$ and $(u_n, v_n) =1$. Note that $y_n=y_0^{(-2)^n}=\sigma(x_n)$ for all $n\geq 0$. Then 
$$ 
x_n-x_m=\frac {  (\alpha\delta-\beta\gamma) (y_n-y_m)}  {(\gamma y_m-\alpha)(\gamma y_n-\alpha) }
=\frac {  (\alpha\delta-\beta\gamma) (u_nv_m-u_mv_n)}  {(\gamma u_m-\alpha v_m)(\gamma u_n-\alpha v_n) }
$$
and in particular  
$$ 
x_{n+1}-x_n 
=\frac {  (\alpha\delta-\beta\gamma) (v_{n}^2v_n-u_nu_n^2)}  {(\gamma u_n-\alpha v_n)(\gamma v_{n}^2-\alpha u_n^2) } 
=\frac {  (\alpha\delta-\beta\gamma) (v_{n}^3- u_n^3)}  {(\gamma u_n-\alpha v_n)(\gamma v_{n}^2-\alpha u_n^2) } 
$$
and
$$ 
x_{n+2}-x_n 
=\frac {  (\alpha\delta-\beta\gamma) (u_{n}^4v_n-u_nv_{n}^4)}  {(\gamma u_n-\alpha v_n)(\gamma u_{n}^4-\alpha v_{n}^4) }
=- \frac {  (\alpha\delta-\beta\gamma)u_nv_n (v_{n}^3- u_n^3)}  {(\gamma u_n-\alpha v_n)(\gamma u_{n}^4-\alpha v_{n}^4) }  .
$$
Let $\eta=\gamma u_n-\alpha v_n$. 
Now $(u_n,v_n)=(1)$ so that $(v_{n}^3- u_n^3,u_n)=(v_{n}^3- u_n^3,v_n)=(1)$, an equality of $R_S$-ideals. Hence 
 $(v_{n}^3- u_n^3, \eta(\gamma v_{n}^2-\alpha u_n^2)) = (v_{n}^3- u_n^3, \eta(\gamma v_{n}^3-\alpha u_n^2v_n)) = (v_{n}^3- u_n^3, \eta(\gamma u_{n}^3-\alpha u_n^2v_n)) 
= (v_{n}^3- u_n^3, \eta^2)$. Similarly $(v_{n}^3- u_n^3, \eta(\gamma u_{n}^4-\alpha v_{n}^4) )=
(v_{n}^3- u_n^3, \eta^2) $. Hence the prime factors in the numerator of
$x_{n+2}-x_n$ are a subset of those in the numerator of
$x_{n+1}-x_n$, and those dividing $u_nv_n$, which are the same as those dividing 
 $u_0v_0$.  We deduce that there can be no doubly primitive prime factor of the numerator of
$x_{n+2}-x_n$ for any $n\geq 2$. This completes the proof for the maps from case (3) of Lemma~5.
  
	Finally, if $\Delta = 2$ or $3$ and $\phi \in \Cal B_{\Delta, d}$, that is cases (4) and (5) of Lemma 5, then Corollary~2 of \S7 shows $x_{n+\Delta} - x_n$ fails to have a doubly primitive prime factor in its numerator when $n$ is large. \qed
	
 \enddemo

\head 10. The density of prime divisors of dynamical sequences \endhead

Given a sequence $\{ x_n\}_{n\geq 0}$ let $\Cal P$ be the set of
primes which divide the numerator of
some non-zero element $x_n$, and $\Cal P(x)$ be the number of
elements of $\Cal P$ up to $x$.
We will prove that for the Fermat numbers, one has $\Cal P(x)\ll
x^{1/2}/\log x$. Although this bound  is small
compared to the total number of primes up to $x$, a simple heuristic
indicates  that the true order of magnitude of $\Cal P(x)$ is probably
some power of $\log\log x$!

Select integers $m$ and $N$ so that $2^m\approx x^{1/2}$ and
$2^N\approx x^{2/3}$.
There are $\ll 2^{m}/m$ prime factors of $F_0F_1\ldots F_{m-1}=2^{2^m}-1$
by the prime number theorem.
Any prime divisor $p_n$ of  $F_n=2^{2^n}+1$ is $\equiv 1 \pmod {2^{n+1}}$.
There are $\leq x/2^n$ integers in this arithmetic progression, and so
$\ll x/2^N$ such primes, in total, with $n\geq N$.
The Brun/Titchmarsh theorem tells us that there are $\ll  x/2^n \log (x/2^n)$
 primes $\equiv 1 \pmod {2^{n+1}}$. up to $x$. If $m\leq n<N$, this is
 $\ll x/2^n \log x$, and so  there are $\ll x/2^m \log x$ such primes
in total. Combining these
 observations yields the claim that $\Cal P(x)\ll x^{1/2}/\log x$.

Presumably if $\{ x_n\}_{n\geq 0}$ is a dynamical sequence, obtained
from a rational function
of degree $d\geq 2$, then it might be possible to prove something like
$\Cal P(x)\ll x^{1-1/d+o(1)}$
(except for certain degenerate cases, such as $x_n$ defined as iterates of
$\phi(t)=(t+p)^2 -p$ for any prime $p$).
We expect that the prime divisors of $x_n$  belong  to an increasingly
sparse sequence as $n$ gets larger, since the $x_n$ are values of the
iterated function
$\phi^{(m)}(t)$ for all $n\geq m$. Some result of this type should be
accessible from a study
of the Galois groups of these extensions. In fact there are several
interesting results in the literature.
First, Odoni [9] showed that for the Euclid numbers $E_n$ (where
$E_0=2$ and $E_{n+1}=E_n^2-E_n+1$)
we have $\Cal P(x)\ll \pi(x)/\log\log\log x$; and then in [10] the
remarkable result that
for ``almost all'' monic $\phi(t)\in\Bbb Z[t]$ of given degree $\geq 2$ and
given height,
$\Cal P(x)=o(\pi(x))$ no matter what the value of $x_0\in\Bbb Z$.
Recently Jones [6]
showed that $\Cal P(x)=o(\pi(x))$, no matter what
the value of $x_0\in\Bbb Z$, for  the polynomials
$\phi(t)=t(t-a)+a, t^2+at-1$ ($a\ne 0,2$), $t^2+a$ ($a\ne -1$),
$t^2-2at+a$ ($a\ne -1,1$),
where $a\in \Bbb Z$.

\head A. Appendix on complex dynamics \endhead

	Here we collect a few results that lie in the realm of complex dynamics on $\Bbb P^1(\Bbb C)$ viewed as a Riemann surface. They are all well-known, and we either point the reader to a proof or give our own if it is brief. 
	
	Let $\phi(t) \in \Bbb C(t)$ be a rational function and let $P \in \Bbb P^1(\Bbb C)$ be a fixed point of $\phi$. The {\sl fixed point multiplicity of $P \in \Bbb C$} is the order of vanishing of $\phi(t) - t$ at $t = P$. If $P = \infty$, we observe that $1 / \phi(1/t)$ has a fixed point at $t=0$, and we define the {\sl fixed point multiplicity of $P = \infty$} to be the order of vanishing of $1/\phi(1/t) - t$ at the origin. Note that $1 / \phi(1/t) = (\sigma^{-1} \circ \phi \circ \sigma)(t)$, where $\sigma(t) = 1/t$ is a fractional linear change of coordinates. 
	
%	Now suppose $\phi(t)$ is a rational function and $P \in \Bbb P^1(\Bbb C)$ is a fixed point of $\phi^{(n)}(t)$ for some positive integer $n$. We write $a_{P}(n)$ for the multiplicity of $P$ as a fixed point of $\phi^{(n)}$. 
	
	Continuing with the notation from the last paragraph, we define the {\sl fixed point multiplier} $\lambda_P$ to be the derivative $\phi'(P)$ if $P \not= \infty$. If $P = \infty$, define the {\sl fixed point multiplier} to be
	$$
		\lambda_\infty = \frac{d}{dt}\left(\frac{1}{\phi(1/t)}\right) \Big|_{t=0}
		=\frac{d}{dt} (\sigma^{-1} \circ \phi \circ \sigma)(t) \Big|_{t=0},
	$$
where as above, $\sigma(t) = 1/t$. 
	
	Next let us suppose that $P \in \Bbb P^1(\Bbb C)$ is a periodic point of (exact) period m; i.e., $\phi^{(m)}(P) =P$ and this relation is false if we replace $m$ by any smaller positive integer. Then $\phi^{(n)}(P) = P$ for any integer $n$ divisible by $m$ since $\phi^{(n)}(P) = (\phi^{(m)} \circ \phi^{(m)} \circ \cdots \circ \phi^{(m)})(P) = P$. Conversely, if $\phi^{(n)}(P) = P$, then $m \mid n$, for if we write $n = mq + r$ for some $0 \leq r < m$, then $\phi^{(r)}(P) = \phi^{(r)}(\phi^{(mq)}(P)) = \phi^{(n)}(P) = P$, which implies $r = 0$ by minimality. 
	
\proclaim{Lemma A.1} Let $\phi(t)$ be a rational function of degree $d \geq 2$ and let $P \in \Bbb P^1(\Bbb C)$ 
	be a fixed point of $\phi(t)$. The fixed point multiplicity of $P$ is greater than one if and only if $\lambda_P = 1$. In particular, if $\lambda_P = 0$, then the fixed point multiplicity of $P$ is exactly~1.
\endproclaim

\demo{Proof} We may assume that $P \not= \infty$ by replacing $\phi(t)$ with $1 / \phi(1/t)$ and replacing $P=\infty$ with $P = 0$ if necessary. Note that by definition this does not affect the fixed point multiplicity or the multiplier. Moreover, we may replace $\phi(t)$ by $\phi(t+P) - P$ in order to assume that $P = 0$. By the chain rule, this does not affect the fixed point multiplier. The fixed point multiplicity is unaffected because $\phi(t+P) - P - t = \phi(t+P) - (t+P)$ has a zero of order $m$ at $t = 0$ if and only if $\phi(t) - t$ has a zero of order $m$ at $t = P$.  

	To prove the first assertion, we expand as a power series about the origin: $\phi(t) = \lambda_0 t + c_2 t^2 + c_3 t^3 + \cdots$. Then $\phi(t) - t = (\lambda_0 - 1)t + c_2t^2 + \cdots$ and the result follows. 
%The second assertion follows from the first.  
\qed

\enddemo

\proclaim{Lemma A.2} A rational function $\phi(t)$ of degree~$d \geq 2$ has exactly $d+1$ fixed points when counted with multiplicity. 
\endproclaim

\demo{Proof} Choose a fractional linear change of coordinates $\sigma(t) = (\alpha t + \beta) / (\gamma t + \delta)$ with $\alpha \delta - \beta \gamma \not= 0$, and consider the new rational function $(\sigma^{-1} \circ \phi \circ \sigma)(t)$. One can choose $\sigma$ so that any three given distinct points of
$\Bbb P^1(\Bbb C)$ are sent to any other three given distinct points; so, for example,
we may assume that $\infty$ is not a fixed point. (As $d \geq 2$ there must exist at least one non-fixed point, say because $\phi(t) - t$ is not the zero function.) Changing coordinates doesn't change the derivatives of $\phi(t)$ at a fixed point (by the chain rule), and hence does not affect the multiplicity of a fixed point. Now $\phi(t) = f(t) / g(t)$ where $\deg f \leq \deg g = d$. The fixed points are exactly the solutions to the equation $\phi(t) = t$, or equivalently, the roots of $f(t) = tg(t)$. This equation has degree $d+1$, and hence exactly $d+1$ roots when counted with multiplicity. \qed	
	
\proclaim{Lemma A.3} Let $\phi(t)$ be a rational function of degree $d \geq 2$ and let $P \in \Bbb P^1(\Bbb C)$ 
	be a fixed point of $\phi(t)$, in which case $P$ is also a fixed point of $\phi^{(2)}(t)$.
	\roster
		\item If the fixed point multiplicity of $P$ for $\phi$ is $\ell > 1$, then the fixed point multiplicity of $P$
			for $\phi^{(2)}$ is $\ell$. 
		\item Suppose the fixed point multiplicity of $P$ is $1$ for $\phi$ and $\ell > 1$ for $\phi^{(2)}$. Then
			$\ell = 2k+1$ for some positive integer $k$. Moreover, there exist $k$ distinct critical points 
			$Q_1, \ldots, Q_k$ for $\phi$  --- i.e., $Q_i \in \Bbb P^1(\Bbb C)$ with $\phi'(Q_i) = 0$ --- 
			such that $\phi^{(m)}(Q_i) \rightarrow P$ as $m \to \infty$. 
	\endroster
\endproclaim

\demo{Proof} We may assume that $P \not= \infty$ by replacing $\phi(t)$ with $1 / \phi(1/t)$ and replacing $P=\infty$ with $P = 0$ if necessary. Note that by definition this does not affect the fixed point multiplicities or multipliers. To prove (1), we expand as a power series about the origin: $\phi(t) -t = c_\ell t^\ell + \cdots$, or $\phi(t) = t + c_\ell t^\ell + \cdots$. Iterating shows $\phi^{(2)}(t) - t =  2c_\ell t^\ell + \cdots$, from which (1) follows.

	The proof of (2) is significantly more difficult. See [1, Lem.~4] or [8, Lem.~10.4 / Lem.~10.11]. \qed
\enddemo
	
\enddemo

\head Appendix B. Kisaka's classification of exceptions to Baker's Theorem \endhead

	Two rational maps $\phi(t), \psi(t) \in \Bbb C(t)$ are said to be {\sl conjugate} if $\phi = \sigma^{-1} \circ \psi \circ \sigma$ for some fractional linear transformation $\sigma(t) = (\alpha t + \beta) / (\gamma t + \delta)$ with $\alpha \delta - \beta \gamma \not= 0$.

	In [7, Thm.~1] Kisaka showed that if $\phi(t) \in \Bbb C(t)$ is a rational map of degree $d \geq 2$ with no periodic point of exact period $\Delta$, then $\phi$ is conjugate to one of the following:
	\roster 
		\item $(\Delta, d) = (2,2)$
			$$
				\psi(t) = \frac{t^2 - t}{at+1} \qquad (a \not= -1);
			$$
			
		\item $(\Delta, d) = (2,3)$
			$$
				\aligned
				\psi(t) &= \frac{t^3 + at^2 - t}{(a^2 - 1)t^2 - 2at + 1} \qquad (a \not= 0); \text{ or} \\
				\psi(t) &= \frac{t^3 - t}{-t^2 + bt + 1} \qquad (b \not= 0); \text{ or} \\
				\psi(t) &= \frac{t^3 + \frac{4}{c}t^2 - t}{-t^2+ct+1} \qquad (c \not= 0, \pm 2i);
				\endaligned
			$$

		\item $(\Delta, d) = (2,4)$
			$$
					\psi(t) = \frac{t^4 - t}{-2t^3+1} \text{ or } \frac{t^4+t^3+t^2-t}{-t^3+t^2-3t+1}	 
						\text{ or }
					\frac{t^4 - 3^{\frac{1}{3}}\cdot t^3 + 3^{\frac{2}{3}}\cdot t^2 - t}
						{-t^3 + 3^{\frac{4}{3}}\cdot t^2
					-5\cdot 3^{-\frac{1}{3}} \cdot t + 1} \text{ or }
				\frac{t^4 + \overline{c}_0 t^3 + \overline{b}_0 t^2 - t}{-t^3+b_0 t^2 + c_0 t + 1},
			$$
		%		where $\displaystyle b_0 = \frac{3 + \sqrt{5}}{2}, \overline{b}_0 = \frac{3 - \sqrt{5}}{2}, 
		%	c_0 = \frac{-5 - \sqrt{5}}{2}, \overline{c}_0 = \frac{-5+\sqrt{5}}{2}$.
		        where $(x - b_0)(x - \overline{b}_0) = x^2 - 3x + 1$ and $(x - c_0)(x-\overline{c}_0) = x^2 + 5x + 5$.

		\item $(\Delta, d) = (3,2)$
			$$
				\psi(t) = \frac{t^2 + \omega t}{\frac{\omega + 5}{\omega - 1} t + 1} \text{ or } 		
				\frac{t^2 + \omega t}{\omega t + 1},
			$$
		where $\omega$ is a primitive third root of unity.

	\endroster

\enddocument